\documentclass[12pt]{amsart} 
\usepackage{amssymb}

\let\oldint\int
\def\int{\oldint\limits}

\makeatletter

\def\stacksub#1#2{_{\scriptstyle #1\atop\scriptstyle #2}}

\let\oldmarginpar\marginpar
\long\def\marginpar#1{\oldmarginpar{\tiny\raggedright#1\par}}
\@mparswitchfalse 
\def\page{\marginpar}

\def\marginpar#1{\ignorespaces}

\usepackage{euscript}

\let\cal\EuScript

\def\cases#1{\left\{\,\vcenter{\normalbaselines\m@th
    \ialign{$##\hfil$&\quad##\hfil\crcr#1\crcr}}\right.}
\def\matrix#1{\null\,\vcenter{\normalbaselines\m@th
    \ialign{\hfil$##$\hfil&&\quad\hfil$##$\hfil\crcr
      \mathstrut\crcr\noalign{\kern-\baselineskip}
      #1\crcr\mathstrut\crcr\noalign{\kern-\baselineskip}}}\,}

\newtheorem{theorem}{Theorem}
\newtheorem{lemma}[theorem]{Lemma}

\theoremstyle{definition}

\newtheorem{remark}{Remark} 
 
\newtheorem{remarks}{Remarks}

\newtheorem{example}{Example} 
\newtheorem{examples}{Examples}

\def\eqlabel#1{\label{eq#1}}
\def\eqref#1{(\ref{eq#1})}
\numberwithin{equation}{section}

\let\<\langle
\let\>\rangle

\newcommand\R{\mathbb R}

\newcommand\I{II}

\def\O{\cal{O}}

\def\span{\operatorname{span}}

\def\span{\operatorname{span}}

\def\tan{\operatorname{tan}}

\def\lim{\operatorname{lim}}

\def\sin{\operatorname{sin}}
\def\cot{\operatorname{cot}}

\def\dist{\operatorname{dist}}

\let\over\@@over
\let\atop\@@atop
\let\above\@@above
\let\overwithdelims\@@overwithdelims
\let\atopwithdelims\@@atopwithdelims
\let\abovewithdelims\@@abovewithdelims

\def\eqalign#1{\null\,\vcenter{\openup\jot
  \ialign{\strut\hfil$\displaystyle{##}$&$\displaystyle{{}##}$\hfil
      \crcr#1\crcr}}\,}

\makeatother

\begin{document}

\title[Singular Integrals Associated to Hypersurfaces: $L^2 $ Theory]
{Singular Integrals Associated to Hypersurfaces:\\$L^2 $ Theory}

\author{Stephen Wainger}
\address{\hskip-\parindent
Stephen Wainger\\
Department of Mathematics\\
University of Wisconsin\\
Madison, WI  53706\\
USA}
\email{wainger@math.wisc.edu}

\author{James Wright}
\address{\hskip-\parindent
James Wright\\
School of Mathematics\\
University of New South Wales\\
Sydney  2052\\
Australia}
\email{jimw@maths.unsw.edu.au}

\author{Sarah Ziesler}
\address{\hskip-\parindent
Sarah Ziesler\\
Department of Mathematics\\
University College Dublin\\
Belfield, Dublin  4\\
Ireland}
\email{zies@ollamh.ucd.ie}

\begin{abstract}
We consider singular integrals associated to a classical
Calder\'on-Zygmund kernel $K$ and a hypersurface given
by the graph of $\varphi(\psi(t))$ where $\varphi$ is an arbitrary
$C^1$ function and $\psi$ is a smooth convex function of finite type.
We give a characterization of those Calder\'on-Zygmund
kernels $K$ and convex functions $\psi$ so that the associated singular integral 
operator is bounded on $L^2$ for all $C^1$ functions $\varphi$. 
\end{abstract}

\thanks{\hskip-\parindent
Wainger was supported in part by an NSF grant. Wright was 
supported in part by an ARC grant. Research at MSRI was supported
in part by NSF grant DMS-9701755.}
\maketitle
\section{Introduction}
\marginpar{\cite{CZ} lacks a title. Authors: fill in?}
The main purpose of this paper is to investigate the $ L^2 $
boundedness of singular integral operators associated to hypersurfaces  
in $ \R^n $, $ n \geq 3 $. Let $ \Gamma(t) $ be a $ C^1 $ mapping from
a neighborhood 
\page{page 2 of original}
of the origin in $ \R^{n - 1} $ into $ \R^n $ with $ \Gamma(0) = 0 $.
For $x$ in $\R^n$ and $ f $  a $ C^1 $ function with compact support in $
\R^n $, we set
$$
H f(x) = \lim\limits_{\epsilon \rightarrow 0} \int_{\epsilon \leq |t| \leq 1}
f(x - \Gamma(t))K(t) \, dt
$$
where $ K(t)$ is a Calder\'on-Zygmund kernel in $ \R^{n - 1} $. That
is, $K$ is smooth $ (C^\infty) $ away from the origin,
$$
\int_{a \leq |t| \leq b} K(t) \, dt = 0
$$
for every $ 0 < a < b $, and 
$$
K(\lambda t) = \lambda^{-n + 1} K(t)
$$
for every $ \lambda >  0 $.
\page{3}

It is known that if $ \Gamma(t) $ is smooth and the vectors $ \{
\frac{\partial^{\alpha} \Gamma}{\partial t^\alpha}(0)\} $, given by
the derivatives of $ \Gamma $ at the origin, $ \span \R^n $, then
$$
\|H f\|_{L^p} \leq A_p \, \|f\|_{L^p}, \quad 1 < p < \infty.
$$
See \cite{St} for this result. Our main interest is studying what
happens when the vectors $
\frac{\partial^{\alpha} \Gamma}{\partial t^{\alpha}}(0) $ do {\em
not} span $ \R^n $. We shall consider surfaces of the form
$$
\Gamma(t) = \big(t, \varphi(\psi(t))\big)
$$
where $ t = (t_1, \ldots, t_{n - 1}) $ and $ \psi(t) $ is
\page{4}
a smooth convex function of finite type with $ \psi(0) = \nabla \psi
(0) = 0 $. (We say that $ \psi(t) $ is of finite type if the graph $
t_n = \psi(t) $ has no lines tangent to infinite order.) If $ \psi(t)
= |t|^2 = t^2_1 + \cdots + t^2_{n - 1} $, then 
$$
\|H f\|_{L^2} \leq A\, \|f\|_{L^2}
$$
for any $ \varphi $. The details of this easy calculation can be found
in \cite{KWWZ}.

The main purpose of this paper is to decide for what convex
\page{5}
functions of finite type $\psi $ and Calder\'on-Zygmund kernels $K$ do
we have
$$
\|Hf\|_{L^2} \leq A \|f\|_{L^2}
$$
for all $ C^1 $ functions $ \varphi $ with $ \varphi(0) = 0 $.
To give the answer to this problem we introduce certain sets which
were considered by Schulz, \cite{Sc}. Let
\marginpar{Is that $s$ in the equation?} 
$$ 
E_{\ell} = \{v \in \R^{n - 1} \mid  \psi(s v) = \cal{O} ( s^{\ell + 1}) \text{
for small } s >  0 \}.
$$
From the convexity of $ \psi $, each $ E_{\ell} $ is a linear subspace
of $ \R^{n - 1} $. Clearly $ E_1 = \R^{n - 1} $, $ E_{\ell + 1}
\subseteq E_{\ell} $ and $ \bigcap E_{\ell} = \{ 0 \} $ (from the
finite type condition). We let $ \ell_0 $ be the smallest value of $
\ell $
\page{6}
such that $ E_{\ell} $ is not all of $ \R^{ n - 1} $. We then have
the following theorem.

\begin{theorem}
\label{1}
If the codimension of $ E_{\ell_0} $ in $ \R^{n - 1} $ is at least $ 2
$, then
$$
\|Hf\|_{L^2} \leq A \, \|f\|_{L^2}
$$
for all $ C^1 $ functions $ \varphi $ with $ \varphi(0) = 0 $.
\marginpar{Silvio recommended ``$ C^1 \varphi $'' $ \rightarrow $ ``$ C^1 $ functions $ \varphi $''}
\end{theorem}

If the codimension of $ E_{\ell_0} $ is $ 1 $, then $H$ is bounded on
$ L^2 $ for all $ C^1 $ functions $ \varphi $ if and only if $K$ satisfies an
additional cancellation condition.
\page{7}
\begin{theorem}
\label{2}
Suppose the codimension of $ E_{\ell_0} $ is $ 1 $, and let $v$ be a
non-zero vector in $ E^\perp_{\ell_0} $. Then
$$
\|Hf\|_{L^2} \leq A\, \|f\|_{L^2}
$$
for all $ C^1 $ functions $ \varphi $ with $ \varphi(0)= 0 $ if and only if $ K(t)
$ satisfies the additional cancellation condition
\begin{equation}
\eqlabel{1.1}
\int\stacksub{v \cdot t \geq 0}{a \leq |t| \leq b} K(t) \,
dt = 0
\end{equation}
for all $ 0 < a < b $.
\end{theorem}
\page{8}
\begin{remarks}

\hspace*{.2in}

\begin{itemize}
\item[(1)]
The positive assertions in Theorems \ref{1} and \ref{2} hold for
the more general operators
\marginpar{Only a mathematician could have guessed $ f(x - \Gamma(+
1)) $ from the handwriting here.}
$$
Hf(x) = \int b(\psi(t)) K(t) f(x - \Gamma(t)) \, dt
$$
for any bounded function $b$, with no change in the proof.
\item[(2)] Theorem \ref{1} is vacuous and Theorem \ref{2} is trivially true when $ n = 2 $,
and so nothing new is being proved for singular integrals along
curves in the plane.
\page{9}
\end{itemize}
\end{remarks}
\begin{examples}

\hspace*{.2in}

\begin{itemize}
\item[(1)] $ \psi(x,y,z) = x^2 + y^2 + z^4 $ is a convex function of
finite type where $ \ell_0 = 2 $ and 
$$ 
E_{\ell_0} = \{  (0,0,z) \mid z \in \R \} 
$$ 
has codimension $ 2 $ and so Theorem 1 applies.
\item[(2)] $ \psi(x,y,z) = x^2 + y^4 + z^4 $ is also a convex function
of finite type where $ \ell_0 = 2 $ but 
$$ 
E_{\ell_0} = \{ (0, y, z) \mid  y, z \in \R \} 
$$ 
has codimension  $ 1 $ and so Theorem 2 applies with $ v = (1, 0,0) $.
\end{itemize}
\end{examples}
\page{10}
\vskip 5pt

Next we turn to examine what happens when the cancellation condition
\eqref{1.1} is not satisfied.

\begin{theorem}
\label{3}
Let $ \bar{\varphi}(s) = \varphi(s^{\ell_0}) $. Assume the codimension
of $ E_{\ell_0} $ is $ 1 $, and the cancellation condition {\rm \eqref{1.1}}
fails. Then if $ \bar{\varphi}(s) $ is convex,
$$
\|Hf\|_{L^2} \leq A \, \|f\|_{L^2}
$$
if and only if 
$$
\bar{\varphi}'(Cs) \geq 2 \bar{\varphi}' (s) 
$$
for some $ C \geq 1 $ and all $ 0 < s \leq 1 $. 
\end{theorem}
\page{11}
\begin{remarks}
\hspace*{.2in}
\begin{itemize}
\item[(1)]
The significance of the power $ \ell_0 $ is
that $\frac{1}{\ell_0} $  is the smallest power $
\alpha $ such that $ [\psi(t)]^{\alpha} $ is a convex function.
\item[(2)]
When $\phi(s)=|s|^2$ and so $n=1$, Theorem 3 was proved in [NVWW].
\end{itemize}
\end{remarks}
\vskip 7pt
If $ E_{\ell_0} = \{ 0 \} $, which means that $ \psi $ is
approximately homogeneous of degree $ \ell_0 $, we obtain $ L^p $
results for $H$ and the corresponding maximal function
$$
Mf(x) = \sup_{0 < h \leq 1} \, \frac{1}{h^{n - 1}} \, \int_{|t| \leq h} |
f(x - \Gamma(t))| \, dt.
$$
We again set $ \bar{\varphi}(s) = \varphi(s^{\ell_0}) $.
\page{12}
\begin{theorem}
\label{4}
Suppose $n \ge 3$, $ E_{\ell_0} = \{ 0 \} $ and $ \bar{\varphi}(s) $ is convex.
Then
$$
\|Hf\|_{L^p} \leq A_p \, \|f\|_{L^p}, \quad 1 < p < \infty ,
$$
and
$$
\|M f\|_{L^p} \leq A_p\, \|f\|_{L^p}, \quad 1 < p \leq \infty.
$$
\end{theorem}
\begin{remark}
It is known that the assertion of Theorem \ref{4} fails in general if
the hypothesis that $ \bar{\varphi} $ is convex is dropped, even if $
\psi(t) = |t|^2 $. See \cite{SWWZ}.
\end{remark}
\page{13}
\vskip 5pt
Finally we make one observation in $ \R^3 $ in the case that $ \psi(t)
$ is not of finite type. Let $ t_0 $ be a point on the curve $ \psi(t)
= 1 $ and $ \ell(t_0 ) $ denote the line tangent to $ \psi(t) = 1 $ at $ t_0
$. Set
$$
E(t_0, \epsilon) = \{ s \in \R^2 \mid  \psi(s) = 1 \text{ and }  \dist (s,
\ell(t_0))\leq \epsilon \}.
$$
\begin{theorem}
\label{5}
Assume $ \psi(t) $ is convex and homogeneous of degree $ 1 $. Then if
$$
\sup\stacksub{t_0}{\psi(t_0) = 1} \int^1_0 |E(t_0,
\epsilon)|\, \frac{d \epsilon}{\epsilon} < \infty, 
$$
$ \|Hf\|_{L^2} \leq A\,\|f\|_{L^2} $ for every $ C^1 $ function $ \varphi $.
\end{theorem}
\page{14}
\vskip 5pt
\begin{example} Consider a smooth convex function $\psi(x,y)$,
homogeneous of degree 1, such that for $|x| << |y|$,
$$\psi(x,y) \ = \ \sqrt{x^2 + y^2} \ {\rm exp}\,\left(- \left(\frac{\sqrt{
x^2 + y^2}}{|x|}\right)^{\alpha}\right) .$$
Clearly $\psi$ is not of finite type and the integrability
condition in Theorem 5 is satisfied exactly when $\alpha < 1$.
\end{example}
\vskip 5pt
In section 2  we will prove Theorems \ref{1} and \ref{2} in the
special cases where $ \psi(x,y,z) = x^2 + y^2 + z^4 $ (for Theorem
\ref{1}) and $ \psi(x,y,z) = x^2 + y^4 + z^4 $ (for Theorem \ref{2}),
where the main direction of the proof is not clouded by intricate
estimates. 
The proof for Theorem
\ref{1} in the general case will
\page{15}
be given in section 3. Theorems \ref{2} and \ref{3} will be proved in
section 4 and sections 5 and 6 contain the proofs of Theorems \ref{4}
and \ref{5} respectively.

Our work is heavily dependent on ideas of Schulz, \cite{Sc}. We would
like to thank A.~Iosevich for bringing the paper \cite{Sc} to our
attention. 
We would also
like to thank Professor A.~Carbery for evaluating a determinant for us.
\page{16}

\section{Special Cases}
In this section we will prove Theorems \ref{1} and \ref{2} in the
special cases $ \psi(x, y, z) = x^2 + y^2 + z^4 $ and $ \psi(x,y,z) =
x^2 + y^4 + z^4 $. We begin with $ \psi(x,y,z) = x^2 + y^2 + z^4 $.
Here no further cancellation condition is required for the
Calder\'on-Zygmund kernel $K$. We need to show
\begin{equation}
\eqlabel{2.1}
\left| \int_{\epsilon \leq x^2 + y^2 + z^2 \leq 1}e^{i \gamma
\varphi(x^2 + y^2 + z^4)} e^{i \eta z} e^{i (\xi_1 x + \xi_2 y)}
K(x,y,z) \, dx\, dy\, dz \right|\leq B 
\end{equation}
uniformly in $ \xi = (\xi_1, \xi_2) , \eta, \gamma $ and $ \epsilon >
0 $. 
\page{17}
Introducing polar coordinates in the $ (x,y) $ integral, the integral
in \eqref{2.1} becomes
\begin{equation}
\eqlabel{2.2}
\int_{\epsilon \leq r^2 + z^2 \leq 1} e^{i \gamma \varphi (r^2 + z^4)}
e^{i \eta z} r \int^{2 \pi}_0 e^{i r(\xi_1 \cos \theta + \xi_2 \sin
\theta )} K(r \cos \theta, r \sin \theta, z) \, d \theta \, dr \, dz.
\end{equation}
We split the integral in \eqref{2.2} as a sum of two integrals $ I_1 +
I_2 $ where the $r$ integration in $ I_1 $ is
\page{18}
restricted to $ r |\xi| \geq 1$ and where the integration in $ I_2 $
is over the complementary range. Using the fact that the
$ \theta $ integral in \eqref{2.2} is the Fourier transform of a
smooth density on the unit circle, we see that
\marginpar{Second fraction here looked like $ r^{1/2}$!}
$$
|I_1| \leq C \int_{r |\xi| \geq 1} \frac{r^2}{(r |\xi|)^{1/2}}
 \int^\infty_{-\infty} \frac{1}{(r^2 + z^2)^{4/2}}+\frac{1}{r(r^3 +|z|^3)
 }\, dz \, dr
= C \int_{r |\xi| \geq 1} \frac{1}{(r|\xi|)^{1/2}} \, \frac{dr}{r} \leq
C.
$$
In $ I_2 $ we replace $ e^{i r (\xi_1 \cos \theta + \xi_2 \sin
\theta)} $ with $ 1 $, creating an error at most a multiple of 
\page{20}
$$
\int_{r |\xi| \leq 1} r^2|\xi| \int^\infty_{-\infty}\frac{1}{(r^2 +
z^2)^{3/2}} \, dz \, dr = C \int_{r|\xi| \leq 1} |\xi| \, dr = C.
$$
Therefore the integral in \eqref{2.2} is
$$
\int^{2 \pi}_0 \!\int \! \int\stacksub{\epsilon \leq r^2 + z^2 \leq
1}{r|\xi| \leq 1} e^{i \gamma \varphi(r^2 + z^4)} e^{i \eta z} K(r
\cos \theta, r \sin \theta, z) r \, dr \, dz \, d \theta \ + \ \O (1).
$$
Furthermore the $ (r, z) $ integration may be further restricted to
the region where $ |z | \leq \delta r^{1/2} $ since integrating $K$ over
the complementary
\page{20}
region is at most
$$
\int \!\int_{\delta r^{1/2} \leq |z| \leq 1} \frac{1}{(r^2 + z^2)^{3/2}}
\, r \, dr \, dz \leq \int_{|z| \leq 1} \frac{1}{|z|^3} \int_{r \leq \left(
\tfrac{1}{\delta} |z|\right)^2} r \, dr \leq C.
$$
With the restriction $ |z| \leq \delta r^{1/2} $ for small $\delta > 0$, we may make the
change of variables $ \lambda = \sqrt{r^2 + z^4} $ (so that $\lambda \sim r$)
in the $r$ integral
to reduce matters to showing that the integral
$$
I = \int^{2 \pi}_{0} \int_{\lambda |\xi| \leq 1} e^{i \gamma
\varphi(\lambda^2)} \lambda \int\stacksub{|z| \leq \delta
\lambda^{1/2}}{\epsilon \leq \lambda^2 + z^2 \leq 1} e^{i \eta z}
K( \sqrt{\lambda^2 - z^4} \cos \theta , \sqrt{\lambda^2 - z^4} \sin
\theta, z) \, dz \, d \lambda \, d \theta
$$
\page{21}
is uniformly bounded in $ \gamma , \eta, \xi $ and $ \epsilon >  0 $.
Replacing $ \sqrt{\lambda^2 - z^4} $ by $ \lambda $ in $I$ creates an
error at most
$$
C \int^1_0 \int_{|z| \leq \lambda^{1/2}} \frac{z^4}{(z^2 +
\lambda^2)^{4/2}} \, dz \, d \lambda \leq C \int^1_0 \lambda^{1/2} \,
d \lambda \leq C
$$
and so
$$
I = \int^{2 \pi}_0 \int_{\lambda|\xi| \leq 1} e^{i \gamma
\varphi(\lambda^2)} \lambda \int\stacksub{|z| \leq \delta
\lambda^{1/2}}{\epsilon \leq \lambda^2 + z^2 \leq 1} e^{i \eta z}
K(\lambda \cos \theta, \lambda \sin \theta, z) \, dz \, d \lambda\, d
\theta \ + \ \O (1).
$$
Next \page{22} we will see that we can replace the oscillatory factor $ e^{i
\eta z} $ with  $ 1 $ in the above integral if we further restrict the $
\lambda $ integration to 
$
\lambda \leq \frac{1}{|\eta|} .
$ 
In fact 
we can integrate by parts in the $z$ integral to
see that the part of the integral where 
$ \lambda |\eta|  \geq 1 $ is at most
$$
C \frac{1}{|\eta|}  \int_{\lambda |\eta| \geq 1} \lambda
\int^\infty_{- \infty}\, \frac{1}{(z^2 + \lambda^2)^{4/2}} \, dz 
\leq C \, \frac{1}{|\eta|} \int_{\lambda |\eta| \geq 1} \,
\frac{1}{\lambda^2} \leq C.
$$
\page{23}
For $ \lambda |\eta| \leq 1$, replacing $
e^{i \eta z} $ by 1 creates an error
at most 
$$
C\, |\eta| \int_{\lambda |\eta| \leq 1} \lambda
\int^\infty_{ - \infty} \frac{|z|}{(\lambda^2 + z^2)^{3/2}} \, dz \, d
\lambda \leq C \, |\eta| \int_{\lambda |\eta| \leq 1} d \lambda \leq C.
$$
Therefore
\page{24}
$$
\eqalign{
I &= \int^{2 \pi}_0 \int\stacksub{\lambda|\xi| \leq 1}{\lambda |\eta| \leq 1} e^{i \gamma \varphi(\lambda^2)} \lambda
\int\stacksub{|z| \leq \delta \lambda^{1/2}}{\epsilon \leq
\lambda^2 + z^2 \leq 1} K ( \lambda \cos \theta, \lambda \sin \theta,
z ) \, dz \, d \lambda \, d \theta \ + \ \O(1)\cr
&= \int^{2 \pi}_0 \int\stacksub{\lambda |\xi| \leq  1}{\lambda |\eta |
\leq 1} e^{i \gamma \varphi (\lambda^2)} \frac{1}{\lambda}
\int\stacksub{|s| \leq \delta \lambda^{-1/2}}{\epsilon \leq \lambda^2
( 1 + s^2) \leq 1} K (\cos \theta, \sin \theta, s)\, ds \, d \lambda
\, d \theta \ + \ \O(1) \cr
&=\int^{2 \pi}_0 \int\stacksub{\lambda|\xi| \leq 1}{\lambda |\eta| \leq
1} e^{i \gamma \varphi (\lambda^2)} \frac{1}{\lambda}
\int_{\epsilon \leq \lambda^2 (1 + s^2) \leq 1}
K(\cos \theta, \sin \theta, s) \, ds \, d
\lambda\, d \theta \ + \ \O(1)\cr
&=- \int\stacksub{\lambda|\xi | \leq 1}{\lambda |\eta| \leq 1} e^{i
\gamma \varphi (\lambda^2)} \frac{1}{\lambda} \int^{2 \pi}_0
\int^\pi\stacksub{0}{\lambda \leq \sin \psi \leq
\frac{\lambda}{\sqrt{\epsilon}}} K(\sin \psi \cos \theta, \sin \psi
\sin \theta, \cos \psi) \sin \psi \, d \psi \, d \theta \, d \lambda +
\O(1).
}
$$
Here we made the change of variables $ z = s \lambda $ followed by $ s
= \cot \psi $ in the 
\page{25}
$z$ integral. Using the fact that 
\marginpar{Is that zero?}
$$
0 \ = \ \int^{2 \pi}_0 \int^{\pi}_0 K(\sin \psi \cos \theta, \sin \psi
\sin \theta, \cos \psi) \sin \psi \, d \psi \, d \theta
$$
we easily see (by splitting the $ \lambda $ integration at $ \lambda =
\sqrt{\epsilon})$  that $I$ is uniformly bounded in $ \gamma , \xi,
\eta $ and $ \epsilon >  0 $. This finishes the proof of Theorem
1 in the case where $\psi(x,y,z) = x^2 + y^2 + z^4$. 
\vskip 5pt

For the example $ \psi(x,y,z) = x^2 + y^4 + z^4 $ we will show that
the integral
\begin{equation}
\eqlabel{2.3}
\int_{\epsilon \leq x^2 + y^2 + z^2 \leq 1} e^{i \gamma \varphi (x^2 +
y^4 + z^4)} e^{i \eta x} e^{i(\xi_1 y + \xi_2 z)} K(x,y,z) \, dx\,
dy\, dz
\end{equation}
is uniformly bounded in $ \gamma, \eta, \xi = (\xi_1, \xi_2) $ and
\page{26}
$ \epsilon >  0 $ under the additional hypothesis that for all $ 0 < a
< b $ 
\begin{equation}
\eqlabel{2.4}
\int\stacksub{a \leq x^2 + y^2 + z^2 \leq b}{x \geq 0} K(x,y, z) \,
dx\, dy\, dz = 0.
\end{equation}
We would like to make the change of variables $ \lambda^2 = x^2 + y^4
+ z^4 $ in the $x$ integral. In order to do this first observe that
the integral in \eqref{2.3} over the region $ \delta|x|^{1/2} \leq
\sqrt{y^2 + z^2} $ is uniformly bounded. In fact
$$
\eqalign{
\int\!\int\!\int_{\delta |x|^{1/2} < \sqrt{y^2 + z^2} \leq 1} & |K(x,y,
z)| \, dx \, dy\, dz \cr
& \leq C\int\!\int_{\sqrt{y^2 + z^2}\leq 1} \,
dy\,dz \int_{\delta|x|^{1/2}\leq \sqrt{y^2 + z^2}}\, \frac{1}{(y^2 +
z^2)^{3/2}} \, dx \cr
& \leq C \int\!\int_{\sqrt{y^2 +  z^2} \leq 1}\, \frac{1}{\sqrt{y^2 + z^2}} \, dy\, dz \leq C.}
$$ 
\page{27}
Hence it suffices to show the uniform boundedness of  
$$
\I \ = \ \int\stacksub{\epsilon \leq x^2 + |\bar{y}| \leq 1}{|\bar{y}| \leq
\delta |x|^{1/2}} e^{i \gamma \varphi (x^2 + y^4 + z^4)} e^{i \eta x}
e^{i \xi \cdot \bar{y}} K(x,\bar{y}) \,dx \, d \bar{y}
$$
where $ \bar{y} = (y,z) $. We write $ \I = \I_+ + \I_- $ where the
integration in $ \I_+ $ is over positive values of $x$. We first
concentrate on $ \I_+ $, making the change of variables
\marginpar{Should $ \I $ be in roman?}
\marginpar{Is this OK?}
$$
\lambda = \sqrt{x^2 + y^4 + z^4}, \quad x = x(\lambda, \bar{y}) =
\sqrt{\lambda^2 - y^4 - z^4}
$$
in the $x$ integral so that $x\sim \lambda$. Then
\page{28}
\marginpar{Ms. illegible at ??}
$$
\I_+ = \int^1_0 e^{i \gamma \varphi (\lambda^2)} \int\stacksub{\epsilon
\leq \lambda^2 + |\bar{y}|^2 \leq 1}{|\bar{y}| \leq \delta
\lambda^{1/2}} e^{i \eta x (\lambda, \bar{y})} e^{i \xi \cdot \bar{y}}
K(x(\lambda, \bar{y}), \bar{y})\, \frac{\partial x}{\partial \lambda} \,
d \bar{y} \, d \lambda \ + \ \O (1).
$$
In order to analyze this integral we make the following simple
observations regarding $ x(\lambda, \bar{y}) $ in the region $
|\bar{y}| \leq \delta \lambda^{1/2}$: 
\hspace*{.2in}
\begin{itemize}
\item[(a)] 
$$ x(\lambda, 0) \ = \ \lambda, \quad\quad  \frac{\partial
x}{\partial \lambda } (\lambda, 0) \ = \ 1 ,$$
\item[(b)] 
$$ \left| \frac{\partial x}{\partial \bar{y}}\right| \ \ \sim \ \
\frac{|\bar{y}|^3}{\lambda}, \quad\quad \left|\frac{\partial x}{\partial \lambda
\partial \bar{y}}\right| \ \ \sim \ \ \frac{|\bar{y}|^3}{\lambda^2} ,  $$
and
\item[(c)] 
$$ \frac{\partial^4 x}{\partial y^4} \ \ \sim \ \ 
\frac{1}{\lambda}, \quad\quad \frac{\partial^4 x}{\partial z^4} \ \sim \ 
\frac{1}{\lambda}. $$
\end{itemize}
\page{29}
Using (a) 
\marginpar{Is this too crowded?}
and (b) we may replace $ \frac{\partial x}{\partial \lambda}
$ by 1 in $ \I_+ $ with an error at most
$$
\eqalign{
C\int^1_0 \,\frac{1}{\lambda^2} \int_{|\bar{y}| \leq \lambda^{1/2}}
\,\frac{|\bar{y}|^4}{\lambda^3 + |\bar{y}|^3} \, d \bar{y} \, d \lambda 
&\leq C \int^1_0 \lambda \int_{|\bar{s}| \leq {\lambda}^{-1/2}}
\,\frac{|\bar{s}|^4}{1 + |\bar{s}|^3} \, d \bar{s}\, d \lambda \cr
& \leq C \int^1_0 \,\frac{1}{\sqrt{\lambda}} \, d \lambda \leq C.
}
$$
Also replacing $ x(\lambda, \bar{y}) $ with $ \lambda $ in the kernel
$K$ creates an error at most
\page{30}
$$
\eqalign{
C \int^1_0 \frac{1}{\lambda} \int_{|\bar{y}| \leq \lambda^{\frac{1}{2}}}\,
\frac{|\bar{y}|^4}{\lambda^4 + |\bar{y}|^4} \, d \bar{y} \, d \lambda
&\leq C \int^1_0 \lambda \int_{|\bar{s}| \leq \lambda^{-1/2}}
\frac{|\bar{s}|^4}{1 + |\bar{s}|^4} \, d \bar{s} \, d \lambda \cr
& \leq C \int^1_0 \, d \lambda = C.
}
$$
Therefore
$$
\I_+ = \int^1_0 e^{i \gamma \varphi(\lambda^2)} \int\stacksub{\epsilon
\leq \lambda^2 + |\bar{y}|^2 \leq 1}{|\bar{y}| \leq \delta
\lambda^{1/2}} e^{i \eta x(\lambda, \bar{y})}
e^{i \xi \cdot \bar{y}} K(\lambda, \bar{y}) \, d
\bar{y} \, d \lambda \ + \ \O(1).
$$
\page{31}
Next we will show that we can replace the oscillation $ e^{i \eta x
(\lambda, \bar{y})} $ with $ e^{i \eta \lambda} $ provided that the $
\lambda $ integration is restricted to where $ \lambda 
|\eta|^{1/3} \leq 1 $. In fact using the fact that
$ \frac{\partial^4 x}{\partial
y^4} \sim \frac{1}{\lambda} $ and Van~der~Corput's lemma (see e.g.,
[St]) in the $y$
integral we see that the part of the integral where $ \lambda 
|\eta|^{1/3} \geq 1 $ is at most
$$
C \frac{1}{|\eta|^{1/4}} \int_{\lambda |\eta|^{1/3} \geq 1} \lambda^{1/4}
\left[ \int \frac{1}{\lambda^4 +
|\bar{y}|^{4}} \, d \bar{y} + \int \frac{1}{\lambda^3 + z^3}\, 
d z \right]\, d \lambda 
\leq C \, \frac{1}{|\eta|^{1/4}} \int_{\lambda |\eta|^{1/3} \geq 1}
\, \frac{\lambda^{1/4}}{\lambda^2} \, d \lambda \leq C.
$$
\page{32}
For the part where $ \lambda |\eta|^{1/3} \leq 1 $ we expect
only to replace $ e^{i \eta x (\lambda, \bar{y})} $ with $ e^{i \eta
\lambda} $ in the region where $ |\bar{y}| \leq
\left(\frac{\lambda}{|\eta|}\right)^{1/4}$ since using (b)
$$
| e^{i \eta x (\lambda , \bar{y})} -e^{i \eta \lambda} | \ \leq \ |\eta|
\, |x (\lambda , \bar{y}) - x (\lambda , 0)| \ \leq \ |\eta| \,
\frac{|\bar{y}|^4}{\lambda}.
$$
In the complementary region, $ \lambda |\eta|^{1/3} \leq 1 $
and $ |\bar{y}| \geq \left(\frac{\lambda}{|\eta|}\right)^{1/4} $ we see
that $K$ is uniformly integrable. In fact
\page{33}
$$
\eqalign{ \int_{\lambda |\eta|^{1/3} \leq 1} \int_{|\bar{y}|
\geq \left(\frac{\lambda}{|\eta|}\right)^{1/4}} |K (\lambda, \bar{y})|
\, d \bar{y} \, d \lambda &\leq \int_{\lambda 
|\eta|^{1/3} \leq 1} \int_{|\bar{y}| \geq
\left(\frac{\lambda}{|\eta|}\right)^{1/4}} \,\frac{1}{\lambda^3 +
|\bar{y}|^3} \, d \bar{y} \, d \lambda \cr & \leq C \int_{\lambda 
|\eta|^{1/3} \leq 1} \, \frac{1}{\lambda} \int_{|\bar{s}| \geq
\left(\frac{\lambda}{|\eta|}\right)^{1/4} \, \frac{1}{\lambda}} \,
\frac{1}{1 + |\bar{s}|^3} \, d \bar{s} \, d \lambda \cr & \leq C
|\eta|^{1/4} \int_{\lambda |\eta|^{1/3} \leq 1} \,
\frac{1}{\lambda^{1/4}} \, d \lambda \leq C.  }
$$
Replacing $ e^{i \eta x (\lambda, \bar{y})} $ with $ e^{i \eta
\lambda} $ in the region $ \lambda |\eta|^{1/3} \leq 1 $ and $
|\bar{y}| \leq \left( \frac{\lambda}{|\eta|}\right)^{1/4} $ creates an
error at most
\page{34}
$$
\eqalign{
&C |\eta| \int_{\lambda |\eta|^{1/3} \leq 1} \, \frac{1}{\lambda}
\int_{|\bar{y}| \leq \left(\frac{\lambda}{|\eta|}\right)^{1/4}}
\, \frac{|\bar{y}|^4}{\lambda^3 + |\bar{y}|^3} \, d \bar{y} \, d \lambda \cr
&\qquad \leq C |\eta| \int_{\lambda |\eta|^{1/3} \leq 1} \lambda^2
\int_{|\bar{s}| = \left(\frac{\lambda}{|\eta|}\right)^{1/4}
\frac{1}{\lambda}} \, \frac{|\bar{s}|^4}{1 + |\bar{s}|^3} \, d
\bar{s} \, d \lambda \cr
& \qquad \leq C \frac{|\eta|}{|\eta|^{3/4}} \int_{\lambda 
|\eta|^{1/3} \le 1} \, \frac{1}{\lambda^{1/4}} \, d \lambda \leq
C.
}
$$
Therefore
$$
\I_+ = \int_{\lambda |\eta|^{1/3} \leq 1} e^{i \gamma
\varphi(\lambda^2)} e^{i \eta \lambda} \int\stacksub{\epsilon \leq \lambda^2 +
|\bar{y}|^2 \leq 1}{|\bar{y}| \leq \delta \lambda^{1/2}} e^{i \xi
\cdot \bar{y}} K(\lambda, \bar{y}) \, d \bar{y} \, d \lambda \ + \ \O(1).
$$
A \page{35} similar but easier argument allows us to replace $ e^{i \xi \cdot
\bar{y}} $ with $ 1 $ if we further restrict the $ \lambda $
integration where $ \lambda |\xi| \leq 1 $. Hence making the
change of variables $ \bar{y} = \lambda \bar{s} $,
$$
\eqalign{
\I_+ &= \int\stacksub{0 \leq \lambda \leq 1}{\lambda \leq
\min(|\xi|^{-1}, |\eta|^{-1/3})}
e^{i \gamma \varphi(\lambda^2)}
e^{i \eta \lambda}
\int\stacksub{\epsilon \leq \lambda^2 + |\bar{y}|^2 \leq 1}{|\bar{y}|
\leq \delta \lambda^{1/2}}
K(\lambda, \bar{y}) \, d \bar{y} \, d \lambda \ + \ \O(1) \cr
&= \int\stacksub{0 \leq \lambda \leq 1}{\lambda \leq \min(|\xi|^{-1},
|\eta|^{-1/3})} 
e^{i \gamma \varphi(\lambda^2)} e^{i \eta \lambda} \, \frac{1}{\lambda}
\, \int_{\epsilon \leq \lambda^2 (1 + |\bar{s}|^2) \leq 1} K(1, \bar{s})
\, d \bar{s} \, d \lambda \ + \ \O (1).
}
$$
\page{36}
Here we used the fact that
$$
\int_{|\bar{s}| \geq \lambda^{-1/2}} |K(1, \bar{s})| \, d \bar{s} \ = \ 
\O(\lambda^{1/2}).
$$
Making a polar change of coordinates $ \bar{s} = ( r \cos \theta, r
\sin \theta ) $ followed by $ r = \tan \psi $, $ 0 \leq \psi < \pi/2
$, in the $ \bar{s} $ integral allows us to write
$$
\eqalign{
\I_+ &= \int\stacksub{0 \leq \lambda \leq 1}{\lambda \leq
\min(|\xi|^{-1}, |\eta|^{-1/3})}
e^{i \gamma \varphi (\lambda^2)} e^{i \eta \lambda} \,
\frac{1}{\lambda} \, \int^{2 \pi}_0 \int\stacksub{0\leq 
\psi \leq \pi/2}{\lambda
\leq \cos \psi \leq \lambda/\sqrt{\epsilon}} K(\cos \psi, \sin
\psi \cos \theta, \sin \psi \sin \theta) \cr\noalign{\vskip - 9pt}
&\hspace*{4in} \sin \psi \, d \psi \, d
\theta \, d \lambda + \O (1) \cr
&= \int\stacksub{\sqrt{\epsilon} \leq \lambda \leq 1}{\lambda \leq
\min(|\xi|^{-1}, |\eta|^{-1/3})}
e^{i \gamma \varphi (\lambda^2)} e^{i \eta \lambda} \,
\frac{1}{\lambda}\,
\int^{2 \pi}_0 \! \int^{\pi/2}_0
K(\cos \psi, \sin \psi \cos \theta, \sin \psi \sin \theta)
\cr\noalign{\vskip - 9pt}
&\hspace{4in} \sin \psi \, d \psi \, d \theta \, d \lambda + \O(1).
}
$$
\page{37}

A similar analysis for $ \I_- $ shows 
$$
\eqalign{
\I_- &= \int\stacksub{\sqrt{\epsilon} \leq \lambda \leq 1}{\lambda \leq
\min(|\xi|^{-1} |\eta|^{-1/3})} e^{i \gamma \varphi (\lambda^2)} e^{-i
\eta \lambda} \, \frac{1}{\lambda}\, \int^{2 \pi}_0 \!\!\int^{\pi}_{\pi/2}
K(\cos \psi, \sin \psi \cos \theta, \sin \psi \sin \theta)
\cr\noalign{\vskip - 6pt}
&\hspace*{4in} \sin \psi \, d \psi \, d \theta\, d \lambda + \O(1).
}
$$
Therefore 
\page{38}
$$
\eqalign{
\I &= \I_+ + \I_- \cr
&= \int\stacksub{\sqrt{\epsilon} \leq \lambda \leq 1}{\lambda \leq
\min(|\xi|^{-1} |\eta|^{-1/3})} 
e^{i \gamma \varphi (\lambda^2)} \sin (\eta \lambda) \,
\frac{1}{\lambda}\, \int^{2 \pi}_0 \!\! \int^\pi_0 K(\cos \psi, \sin
\psi \cos \theta, \sin \psi \sin \theta) \cr\noalign{\vskip - 6pt}
&\hspace*{4in} \sin \psi \, d x \, d \theta \, d \lambda + \O(1)\cr
&= \ 0 \ + \ \O(1)
}
$$
by \eqref{2.4}. Note that when the additional cancellation condition
for $K$ is not satisfied, we are left with a truncated Hilbert
transform along the curve $ (\lambda, \varphi (\lambda^2)) $ and so we
might expect to be able to use the analysis in \cite{NVWW} when $ \varphi 
(\lambda^2) $ is convex.
\page{38}

\section{Proof of Theorem \ref{1}}
According to Schulz \cite{Sc}, after a rotation of coordinates, we may
write 
\begin{equation}
\eqlabel{3.1}
\psi(t) = P(t) + R(t)
\end{equation}
where
$$
P(t) = \sum^r_{j = 1} a_j t_j^{\ell_0} + \sum^{n - 1}_{j = r + 1} a_j
t_j^{m_j} + P_1(t) .
$$
$ P(t) $ is a convex polynomial, $ P(t) >  0 $ for $ t \neq 0 $, $ a_j
>  0 $ for $ 1 \leq j \leq n - 1 $, $ \ell_0 < m_j $ for $ r + 1 \leq
j \leq n - 1 $, $ P_1(t) $ has no pure powers of $t$, and if $ A
t_1^{\alpha_1} \ldots t^{\alpha_{n - 1}}_{n - 1} $ is a monomial
\page{40}
of $ P_1(t) $,
$$
\frac{1}{\ell_0} \sum^r_{j = 1} \alpha_j + \sum^{n - 1}_{j = r + 1}
\, \frac{\alpha_j}{m_j} = 1 .
$$
$ R(t) $ is smooth and if $ At_1^{a_1} \ldots t^{a_{n -
1}}_{n - 1} $ is a term in the Taylor expansion of $ R(t) $ 
$$
\frac{1}{\ell_0} \sum^r_{j = 1} a_j + \sum^{n - 1}_{j = r + 1}
\frac{a_j}{m_j} >  1 .
$$

To prove our theorems, we may assume $ \psi(t) $ has the form
\eqref{3.1}. The hypothesis of Theorem 1 asserts that $ r \geq 2 $.
Let $ H(t) $ be the
\page{41}
part of $ P(t) $ which is homogeneous of degree $\ell_0 $. Then $ H(t)
$ is a function of only $ t_1, \ldots, t_r $. In fact if
$$
A \, t^{\alpha_1}_1 \ldots t^{\alpha_r}_r t^{\alpha_{r + 1}}_{r + 1}
\ldots t^{\alpha_{n - 2}}_{n - 2} t^{\ell_0 - (\alpha_1 + \cdots +
\alpha_{n - 2})}_{n - 1}
$$
were a monomial of $H$, then
$$
\frac{1}{\ell_0} \sum^r_{j = 1} \alpha_j + \sum^{n - 2}_{j = r + 1}
\frac{\alpha_j}{m_j} + \frac{\ell_0 - (\alpha_1 + \cdots + \alpha_{n -
2})}{m_{n - 1}} = 1.
$$
This identity would clearly hold if $ m_{r + 1} = \ldots m_{n - 1} =
\ell_0 $, so it could not hold if one of the $m$'s were bigger than $
\ell_0 $. 
\page{42}
Similarly every monomial of $ P(t) $ which depends only on $ t_1,
\ldots, t_r $ belongs to $H$. So 
$$ 
H(t) = H(t_1, \ldots, t_r) = P( t_1, \ldots, t_r, 0, \ldots, 0) 
$$ 
is convex and positive if some $ t_j $ is nonzero.

We write $ y = (t_1, \ldots, t_r) $ in $ \R^{r} $ and $ x = (t_{r + 1},
\ldots, t_{n - 1}) $ in $ \R^{n - 1 - r} $. We shall suppose $ n - r -
1 \geq 1 $, otherwise the proof is similar but simpler. We then write
$$
P (x,y) = H(y) + P_2 (x,y).
$$
To prove Theorem \ref{1}, we must show
\page{43}
for $ \xi $ in $ \R^r $ and $ \eta $ in $ \R^{n - 1 - r}$,
\begin{equation}
\eqlabel{3.2}
\left| \int_{\epsilon \leq |x|^2 + |y|^2 \leq 1}
e^{i \gamma \varphi(\psi (x,y))} e^{i \eta \cdot x} e^{i \xi \cdot y}
K(y,x) \, dy \, dx \right| \leq B
\end{equation}
uniformly in $ \xi, \eta, \gamma $ and $ \epsilon >  0 $.

We begin by introducing polar coordinates in the $y$ variables. That
is, we write $ y = s \omega $ where $s$ goes from $ 0 $ to $ 1 $ and $
\omega $ runs over the surface $ H(w) = 1 $. The integral in
\eqref{3.2} becomes
\page{44}
\begin{equation}
\eqlabel{3.3}
\int_{H(\omega) = 1}
\int_{\epsilon \leq |x|^2 + s^2 |\omega|^2 \leq 1}
e^{i \gamma \varphi(\psi(x , s \omega))} e^{i( s \xi \cdot \omega +
\eta \cdot x)}
K( s \omega , x) s^{r - 1} h (\omega) \, d s \, d x \, d \omega
\end{equation}
where $h$ is a smooth function. We let $m$ be the smallest of the values
among the $ m_j $'s, $ r + 1 \leq j \leq n - 1 $, and choose $
\sigma >  0 $ so that $ \frac{\ell_0}{m} + \sigma < 1 $. We may
restrict the integration in \eqref{3.3} to $ |x| \leq s^{\ell_0/m +
\sigma} $ since integrating $K$ over the complementary region is at
most
\page{45}
$$
\eqalign{
C \int\!\int_{s^{\frac{\ell_0}{m} + \sigma} \leq |x| \leq 1} \,
\frac{1}{s^{n - 1} + |x|^{n - 1}} s^{r - 1}\, ds \, dx
& \leq C \int_{|x| \leq 1} \, \frac{1}{|x|^{n - 1}} \,\int_{s \leq
|x|^\lambda} \, s^{r - 1} \, ds \, dx \cr
& \leq C \int_{|x| \leq 1} \, \frac{1}{|x|^{n - r \lambda - 1}} \, dx
\leq C
}
$$
since 
$$
 \lambda = \frac{1}{\frac{\ell_0}{m} + \sigma} >  1 \quad \text{ and }
\quad x \in \R^{n - r - 1}.
$$

Furthermore in the region $ |x| \leq s^{\frac{\ell_0}{m} + \sigma} $,
\begin{equation}
\eqlabel{3.4}
\psi(x, s \omega) = s^{\ell_0} \ + \ \O(s^{\ell_0 + \sigma})
\end{equation}
and
\begin{equation}
\eqlabel{3.5}
\frac{\partial \psi}{\partial s} (x , s \omega) = \ell_0 s^{\ell_0 -
1} \ + \ \O (s^{\ell_0 - 1 + \sigma}).
\end{equation}

\noindent In fact
$$
\psi(x, s \omega) = s^{\ell_0} + P_2 (x, s \omega) + R(x, s \omega)
$$
and every monomial in $ P_2 $ or any
\page{46}
monomial in $R$ of the form $ x^\alpha (s \omega)^\beta $ with $ |\alpha| >  0 $
has the bound
$$
|x^\alpha (s \omega)^\beta | \leq s^{\left(\frac{\ell_0}{m} + \sigma\right)
|\alpha| + |\beta|} \leq s^{\ell_0 + \sigma}.
$$
Also any monomial in $R$ of the form $ (s \omega)^\beta $ is $
\O(s^{\ell_0 + 1}) $. Therefore we may make the change of variables
\begin{equation}
\eqlabel{3.6}
\lambda^{\ell_0} = \psi(x , s \omega)
\end{equation}
in $s$ for fixed $x$ and $ \omega $ and write \eqref{3.3} as 
\begin{equation}
\eqlabel{3.7}
\int_{H(\omega) = 1} h(\omega) \int^1_0 e^{i \gamma \varphi(\lambda^{\ell_0})}
\int\stacksub{\epsilon \leq |x|^2 + s^2 |\omega|^2 \leq 1}{|x| \leq
s^{\frac{\ell_0}{m} + \sigma}}
e^{i \eta \cdot x}
e^{i s\xi \cdot \omega} K(s \omega, x) s^{r - 1} \frac{\partial s}{\partial
\lambda} \, dx \, d \lambda \, d \omega + \O(1)
\end{equation}
where $ s = s (\lambda, \omega, x) $. 
\page{47}
From \eqref{3.4} and \eqref{3.5} we have, for some $\epsilon > 0$,
$$
s = \lambda + \O(\lambda^{1 + \epsilon}) \quad \text{ and } \quad
\frac{\partial s}{\partial \lambda} = 1 + \O (\lambda^\epsilon).
$$
Also we have
\begin{equation}
\eqlabel{3.8}
\left| \frac{\partial s}{\partial x} \right| \leq C \lambda^{1 -
\frac{\ell_0}{m}}
\end{equation}
which follows by differentiating \eqref{3.6} with respect to $x$,
giving
$$
0 = \frac{\partial \psi}{\partial x} + \frac{\partial \psi}{\partial
s} \, \frac{\partial s}{\partial x}.
$$
\eqref{3.5} implies $ \frac{\partial \psi}{\partial s} \sim \lambda^{\ell_0
- 1} $ and a similar argument which established \eqref{3.4} and \eqref{3.5} 
shows $
\frac{\partial \psi}{\partial x} = \O (\lambda^{\ell_0 -
\frac{\ell_0}{m}}) $ and thus  \eqref{3.8}.

Arguing as in section 2, using
\page{48}
the above estimates on the derivatives of $ s(\lambda, \omega, x) $
and $s$ itself, shows that we may replace $ s(\lambda, \omega, x) $ by
$ \lambda $ (except in the oscillation $ e^{i s \xi \cdot \omega}$)
and $ \frac{\partial s}{\partial \lambda} $ by $ 1 $ with a uniformly
bounded error. Thus the integral in \eqref{3.7} is
\begin{equation}
\eqlabel{3.9}
\int^1_0 e^{i \gamma \varphi(\lambda^{\ell_0})} \lambda^{r - 1}
\int_{|x| \leq \lambda^{\frac{\ell_0}{m} + \sigma}} e^{i \eta \cdot x}
\int\stacksub{H(\omega) = 1}{\epsilon \leq |x|^2 + \lambda^2|\omega|^2
\leq 1} 
e^{i s (\lambda, \omega, x) \xi \cdot \omega}
K(\lambda \omega, x) h(\omega) \, d \omega \, d x \, d \lambda +
\O(1).
\end{equation}
Consider first the contribution to \eqref{3.9} from those values of
$ \lambda $ where $ \lambda|\xi| \geq 1 $. Since $ H(\omega) = 1 $ is
of
\page{49}
finite type we may for each $ \omega_0 $ on $ H(\omega) = 1 $
parametrize $ H(\omega) = 1 $ in a neighborhood of $ \omega_0 $ as
$$
\omega_0 + (\tau_1, \ldots, \tau_{r - 1}, g (\tau_1, \ldots, \tau_{r -
1}))
$$
where $ g(0) = 0 $, $ \nabla g (0) = 0 $, and for some $ j_0 \geq 2 $,
$$ 
\frac{\partial^j g}{\partial \tau_1^j} (0) = 0 
$$ 
for $ 1 \leq j \leq
j_0 - 1 $ and 
$$ 
\frac{\partial^{j_0} g }{\partial \tau_1^{j_0}} (0)
\neq 0.
$$
It follows that we may assume 
$$ 
\frac{\partial^{j_0} g}{\partial \tau^{j_0}_1} \neq 0 
$$ 
for all $ \tau $ in a neighborhood of $0 $. Therefore since
$s(\lambda, \omega, x) \sim \lambda$,
$$
\left|\int\stacksub{|\omega - \omega_0 | \leq \delta}{H(\omega) = 1}
e^{i s(\lambda, \omega, x) \xi \cdot \omega} h(\omega) \, dw \right|
\leq C \frac{1}{(\lambda|\xi|)^\delta}
$$
for \page{50}some positive $ \delta $ by Van~der~Corput's lemma .
Integrating by parts now shows
that the contribution to the integral in \eqref{3.9} from those $
\lambda $ where $ \lambda |\xi| \geq 1 $ is at most
$$
C \int_{\lambda|\xi| \geq 1} \lambda^{r - 1} \frac{\lambda}{(\lambda
|\xi|)^\delta} \int_{\R^{n - r - 1}} \frac{1}{\lambda^n + |x|^n} \, d
x \, d\lambda 
\leq C \int_{\lambda |\xi| \geq 1} \, \frac{1}{(\lambda |\xi|)^\delta}
\, \frac{d\lambda}{\lambda} \leq C.
$$

Thus the proof of Theorem \ref{1} reduces to showing that the integral
\page{51}
\marginpar{Check subscripts here}
$$
I = \int_{\lambda |\xi| \leq 1}
e^{i \gamma \varphi(\lambda^{\ell_0})} \lambda^{r - 1}
\int_{H(\omega)=1}\!
\int\stacksub{\epsilon \leq |x|^2 + \lambda^2
|\omega|^2 \leq 1}{|x| \leq \lambda^{\frac{\ell_0}{m} + \sigma}}
e^{i \eta \cdot x}
K(\lambda\omega, x) h(\omega) \, dx \, d \omega\, d \lambda
$$
is uniformly bounded in $ \gamma, \xi, \eta $ and $ \epsilon >  0 $.
Putting $ x = \lambda z $ makes
$$
I = \int_{\lambda |\xi| \leq 1}
e^{i \gamma \varphi(\lambda^{\ell_0})} \, \frac{1}{\lambda}\,
\int_{H(\omega)=1} \!\int\stacksub{\epsilon \leq \lambda^2 (|z|^2 +
|\omega|^2 ) \leq 1}{|z| \leq \lambda^{\frac{\ell_0}{m} + \sigma - 1}}
e^{i \lambda \eta \cdot z} K(\omega, z) h(\omega) \, dz \, d \omega \,
d \lambda
$$
and using the fact
$$
\int_{\left(\frac{C}{\lambda}\right)^\delta \leq |z|} |K (\omega, z)|
\, dz \ = \ \O (\frac{\lambda}{C})^{r \delta}
$$
three times, first 
with $ \delta = 1 $ and $ C = \sqrt{\epsilon} $, then with $\delta = 1$
and $ C = 1 $, and finally with
\page{52}
$ \delta = 1 - \left(\left(\frac{\ell_0}{m}\right) + \sigma\right) $  
and $ C = 1 $, we
see that
$$
I = 
\int_{A \sqrt{\epsilon} \leq \lambda  \leq \frac{1}{|\xi|}}
e^{i \gamma \varphi(\lambda^\ell_0)}\, \frac{1}{\lambda}
\int_{H(\omega)= 1}\!\int
e^{i \lambda \eta \cdot z}
K(\omega, z) h(\omega) \, dz \, d \omega \, d \lambda \ + \ \O(1)
$$
if $A$ is chosen large enough. An integration by parts in the $z$
integral shows that the part of the integral 
where $ \lambda |\eta|
\geq 1 $ is at most (up to boundary terms)
$$
C \frac{1}{|\eta|} \int_{\lambda |\eta| \geq 1} \, \frac{1}{\lambda^2}
\, 
\int_{\R^{n - r - 1}} \sup_\omega |\nabla K (\omega, z)| \, dz \, d
\lambda
\leq C \frac{1}{|\eta|} \int_{\lambda |\eta| \geq 1} \,
\frac{1}{\lambda^2} \, 
\int_{|z| \geq 1} \, \frac{1}{|z|^n} \, dz \, d \lambda \leq C,
$$
\page{53}
and so
$$
I = \int\stacksub{A \sqrt{\epsilon} \leq \lambda \leq 1}{\lambda \leq
\min(\frac{1}{|\xi|}, \frac{1}{|\eta|})}
e^{i \gamma \varphi(\lambda^{\ell_0})} \, \frac{1}{\lambda} \, 
\int{H(\omega) = 1}\!\int e^{i \lambda \eta \cdot z} K(\omega, z)
h(\omega) \, dz d \omega \, d \lambda \ + \ \O(1).
$$
The boundary terms are handled similarly. 
Replacing $ e^{i \lambda \eta \cdot z} $ by $ 1 $ creates an error at
most
$$
C |\eta| \int_{\lambda |\eta| \leq 1} \int_{H(\omega) = 1}\!\int
\frac{|z|}{|\omega|^{n - 1} + |z|^{n - 1}} \, dz \,d \omega\, d
\lambda \leq C
$$
since $ r \geq 2 $. Therefore
\page{54}
$$
I = \int\stacksub{A \sqrt{\epsilon} \leq \lambda \leq 1}{\lambda \leq
\min(\frac{1}{|\xi|}, \frac{1}{|\eta|})}
e^{i \gamma \varphi(\lambda^{\ell_0})}\, \frac{1}{\lambda}\,
\int_{H(\omega) = 1}\!\int
K(\omega, z) h(\omega) \, dz d \omega\, d \lambda \ + \ \O(1),
$$
and so it suffices to show
\begin{equation}
\eqlabel{3.10}
\int_{H(\omega) = 1} \int_{\R^{n - r - 1}}
K(\omega, z) h(\omega) \, d z \, d \omega = 0.
\end{equation}
Now for any $ \delta >  0 $
\page{55}
$$
\eqalign{
0 & = \int_{1 - \delta \leq |y|^2 + |x|^2 \leq 1}
K(y, x) \, dy\, dx \cr
&= \int_{H(\omega) = 1} \int_{1 - \delta  \leq \lambda^2 |\omega|^2 +
|x|^2 \leq 1}
\lambda^{r - 1} K(\lambda \omega, x) h(\omega) \, d \lambda \, d
\omega \, dx \cr
&= \int_{H(\omega) = 1}\int_{1 - \delta \leq \lambda^2 (|\omega|^2 +
|z|^2) \leq 1}
K(\omega, z) h(\omega) \, \frac{d \lambda}{\lambda}\, d \omega\, dz \cr
&= \int_{H(\omega)=1} h(\omega) \int_{\R^{n - r - 1}} K(\omega, z)
\int_{\frac{1 - \delta}{|\omega|^2 + |z|^2} \leq \lambda^2 \leq
\frac{1}{|\omega|^2 + |z|^2}} \, \frac{d \lambda}{\lambda} \, dz \, d
\omega.
}
$$
Dividing by $ \delta $ and  letting $ \delta \rightarrow 0 $ gives
\eqref{3.10} and this finishes the proof of Theorem \ref{1}.
\qed

\section{The proof of Theorems \ref{2} and \ref{3}.}
We may again assume $ \psi(t) $ is of the general form \eqref{3.1},
where now $ r = 1 $. The cancellation condition \eqref{1.1} now becomes
$$
\int_{\Sigma^+} K(t)\, d \sigma(t) = 0
$$
where $ \Sigma^+ = \{ t \in \R^{n - 1} \mid |t| = 1, t_1 >  0 \} $ is
the ``upper'' hemisphere of $ S^{n - 2} $. It will be convenient to
let $ t_1 $ be denoted by $y$ and $ (t_2, \ldots, t_{n - 1}) = x \in
\R^{n - 2}  $. We are then concerned with the uniform boundedness of
\page{57}
$$
\eqalign{
\int\int e^{i \gamma \varphi(\psi(y,x))} e^{i \xi \cdot x} e^{i \eta y}
K(y,x) \, dy \, dx &= \int\int_{y >  0} + \int\int_{y < 0} \cr
& = I \ + \ \I.
}
$$
Theorems \ref{2} and \ref{3} will then follow if we can prove for some
$b$, $ 0 < b < 1 $,
\begin{equation}
\eqlabel{4.1}
I = \int\stacksub{0 \leq \lambda \leq 1}{\lambda \leq
\min(\frac{1}{|\xi|}, \frac{1}{|\eta|^b})}
e^{i \gamma \varphi(\lambda^{\ell_0})} e^{i \eta q (\lambda)} \,
\frac{d \lambda}{\lambda} \int_{\Sigma^+} K(\omega) d \sigma(\omega) \ +
\ \O(1),
\end{equation}
\page{58}
\begin{equation}
\eqlabel{4.2}
\I = - \int\stacksub{o \leq \lambda \leq 1}{\lambda \leq
\min(\frac{1}{|\xi|}, \frac{1}{|\eta|^b})}
e^{i \gamma \varphi (\lambda^{\ell_0})}
e^{-i \eta q (\lambda)} \, \frac{d \lambda}{\lambda} \, 
\int_{\sum^+} K(\omega) d \sigma(\omega) \ + \ \O(1),
\end{equation}
and for $ \bar{\varphi}(\lambda) = \varphi(\lambda^{\ell_0}) $ convex,
\begin{equation}
\eqlabel{4.3}
\int\stacksub{0 \leq \lambda \leq 1}{\lambda \leq \min
(\frac{1}{|\xi|}, \frac{1}{|\eta|^b})}
e^{i \gamma \bar{\varphi} (\lambda)} \sin (\eta q (\lambda)) \,
\frac{d \lambda}{\lambda}
\end{equation}
is uniformly bounded in $\gamma, \eta$ and $\xi$ if and only if
$$
\bar{\varphi}' (C \lambda) \geq 2 \bar{\varphi}'(\lambda)
$$
for some $ C \geq 1 $ and $ 0 < \lambda \leq 1 $. Here
$q(\lambda) = \lambda + \O(\lambda^{1 + \epsilon}) $ and $ q '(\lambda) = 1 +
\O(\lambda^\epsilon) $.

We begin with the proof of
\page{59}
\eqref{4.1}. It will convenient to write \eqref{3.1} in the form
\begin{equation}
\eqlabel{4.4}
\psi (y,x) = A y^{\ell_0} + \sum^{n - 2}_{j = 1} a_j x_j^{m_j} +
\sum^{n - 2}_{j = 1} b_j x_j^{\alpha_j} y^{\beta_j} + P_2 (y,x) + R
(y,x)
\end{equation}
where $ {\ell}_0 < m_j $, $ 1 \leq j \leq n - 2 $, $ A >  0 $, $ a_j >  0
$, $ b_j \neq 0 $, $ 1 \leq j \leq n - 2 $ and each monomial of $ P_2 (y,x)
$ has the form $ x^\alpha_j y^\beta $ with $ \alpha >  \alpha_j $, or
contains powers of at least two different $ x_j $'s. Let $ m = \min_{1
\leq j \leq n - 2} (m_j) $ and choose $ \sigma >  0 $ such that $
\frac{\ell_0}{m} + \sigma < 1 $. Again
\page{60}
$$
\int_{|x| \geq y^{\frac{\ell_0}{m} + \sigma}}
|K  (x,y)| \, dx\, dy \ = \ \O(1)
$$
and so it suffices to study $ I $ in the region $ |x| \leq
y^{\frac{\ell_0}{m} + \sigma} $. In this region we wish to make a
change of variables
\begin{equation}
\eqlabel{4.5}
\lambda^{\ell_0} = \psi (y,x) 
\end{equation}
in the $y$ integral. As in section 3, $ |x| \leq y^{\frac{\ell_0}{m} +
\sigma} $ implies that $ y = y(x, \lambda) $ defined implicitly by \eqref{4.5}
satisfies 
\begin{equation}
\eqlabel{4.6}
y(x, \lambda) = A^{- \frac{1}{\ell_0}} \lambda \ + \ \O (\lambda^{1 +
\sigma}),
\end{equation}
\marginpar{Are these $ \O $'s or $ \theta $'s?}
\begin{equation}
\eqlabel{4.7}
\frac{\partial y}{\partial \lambda} = A^{- \frac{1}{\ell_0}} \ +
\ \O(\lambda^\sigma),
\end{equation}
and
\begin{equation}
\eqlabel{4.8}
\left| \frac{\partial y}{\partial x}\right| \leq C \lambda^{1 -
\frac{\ell_0}{m}}.
\end{equation}
\page{61}
Therefore, as before, making the change of variables \eqref{4.5},
shows
\begin{equation}
\eqlabel{4.9}
I = \int_{0 \leq \lambda \leq 1}
e^{i \gamma \varphi (\lambda^{\ell_0})}
\int_{|x| \leq y^{\frac{\ell_0}{m} + \sigma}}
e^{i x \cdot \xi} e^{i \eta y (x, \lambda) } K(x, \lambda) \, d x \, d
\lambda \ + \ \O(1).
\end{equation}
To study \eqref{4.9}, it is necessary to have information on the
derivatives of $y$ with respect to the $x$ variables.
\begin{lemma}
\label{1.l}
Suppose $ |x| \leq \lambda^{\frac{\ell_0}{m} + \sigma} $.
\begin{itemize}
\item[(1)] For $ \delta > 0 $ small,
$$
\left|\frac{\partial^k y}{\partial x_j^k}(x,\lambda)\right| \leq \delta
\lambda^{1 - k \frac{\ell_0}{m_j}}, \ \quad 1 \leq k \leq \alpha_j - 1, \ \ 
1 \leq j
\leq n - 2 .
$$
\page{62}
\item[(2)] 
$$ \frac{\partial^{\alpha_j} y}{\partial x_j^{\alpha_j}}(x,\lambda)
\ \sim \ \lambda^{1 - \alpha_j \frac{\ell_0}{m_j}}, \quad 1 \leq j \leq n
- 2 .$$
\item[(3)] For every $ \beta = ( \beta_1, \ldots, \beta_{n - 2} ) $
with $ 0 \leq \beta_j \leq \alpha_j , \quad 1 \leq j \leq n - 2 $,
$$
\left|\frac{\partial^{\beta}y}{\partial x^{\beta}}(x,\lambda)
\right| \ \leq \ C \lambda^{1 - \ell_0
\sum^{n - 2}_{j = 1} \frac{\beta_j}{m_j}}.
$$
\item[(4)] For every $ \beta = (\beta_1, \ldots, \beta_{n - 2}) $ with
$ 0 \leq \beta_j \leq \alpha_j - 1 $, $ 1 \leq j \leq n - 2 $, either
$$
\frac{\partial^{\beta} y}{\partial x^{\beta}}
(0, \lambda) \ \sim \ \lambda^{p_\beta}
$$
for some $ p_\beta >  -|\beta| $, or
$$
\frac{\partial^{\beta}y}{\partial x^{\beta}}
(0, \lambda) \ = \ \O(\lambda^N)
$$
\page{63}
for every $N$. 
\end{itemize}
\end{lemma}

\begin{proof}[Proof of lemma]
For $M$ large, write
$$
\psi (y,x) = A y^{\ell_0} + \sum^{n - 2}_{j = 1} b_j x_j^{\alpha_j}
y^{\beta_j} + \sum_{u, \beta} c_{u,\beta} x^u y^\beta + \O(|x|^M) + \O(y^M)
$$
where $ A >  0 $, $ b_1, \ldots, b_{n - 2} \neq 0 $, $
\frac{\alpha_j}{m_j} + \frac{\beta_j}{\ell_0} = 1 $ for $ 1 \leq j
\leq n - 2 $, each $ x^u y^\beta $ satisfies $ \sum^{n - 2}_{j = 1}
\frac{u_j}{m_j} + \frac{\beta}{\ell_0} \geq 1 $ and if $ u = (0,
\ldots, u_j, \ldots , 0) $, then $ u_j >  \alpha_j $. To prove (1), we
will show inductively that in the larger region, $ |x_j| \leq \epsilon
\lambda^{\frac{\ell_0}{m_j}} $, $ 1 \leq j \leq n - 2 $,
\page{64}
\begin{equation}
\eqlabel{4.10}
\left| \frac{\partial^k y}{\partial x_j^k}\right| \leq \delta
\lambda^{1 - k \frac{\ell_0}{m_j}}, \quad 1 \leq k \leq \alpha_j - 1,
\end{equation}
provided $ \epsilon = \epsilon (\delta) >  0 $ is small enough. We
first prove \eqref{4.10} for $ k = 1 $ (so $ \alpha_j \geq 2 $ or
there is nothing to prove). If we differentiate \eqref{4.5} with
respect to $ x_j $, noting $ y \sim \lambda $ from \eqref{4.6}, we
obtain
$$
0 \ = \ C_1 \frac{\partial y}{\partial x_j} + C_2
$$
where
$$
C_1 \ = \  A \ell_0 y^{\ell_0 - 1} + \O(\lambda^{\ell_0 + 1}) + E
$$
and  $E$ is a finite sum of terms of the form $ x^u y^{\beta - 1} $
where $ u = (u_1, \ldots, u_{n - 2}) \neq 0 $ and $ \sum^{n - 2}_{j =
1}  \frac{u_j}{m_j} + \frac{\beta}{\ell_0} \geq 1 $. Hence for
\page{65}
$ \epsilon >  0 $ small enough, $ |x^u y^{\beta - 1}| \leq \delta
\lambda^{\ell_0 - 1} $ and therefore
$$
C_1 \sim \lambda^{\ell_0 - 1}.
$$
$ C_2 $ is $ \O(\lambda^{\ell_0}) $ plus a finite sum of terms of the form
$ x^{-1}_j x^u y^\beta $ with $ u_j \geq 2 $ and $ \frac{|u|}{m} +
\frac{\beta}{\ell_0} \geq 1 $. Thus for $ \epsilon >  0 $ small
enough,
$$
\eqalign{
|x^{-1}_j x^u y^\beta| &\leq \delta \lambda^{|u| \frac{\ell_0}{m} +
\beta} \lambda^{- \frac{\ell_0}{m_j}} \cr
& \leq \delta \lambda^{\ell_0 - \frac{\ell_0}{m_j}}
}
$$
and so $ C_2 = \O(\delta \lambda^{\ell_0 - \frac{\ell_0}{m_j}}) $
which proves \eqref{4.10} with $ k = 1 $.

Next we assume \eqref{4.10} for
\page{66}
$ k \leq k_0 - 1 $ where $ k_0 \leq \alpha_j - 1 $, and prove
\eqref{4.10} for $ k = k_0 $. Differentiating \eqref{4.5} $ k_0 $
times with respect to $ x_j $, we again obtain
$$
0 \ = \ D_1 \, \frac{\partial^{{k_0}}y}{\partial x_j^{k_0}} \ + \ D_2
$$
where as before $ D_1 \sim \lambda^{\ell_0 - 1} $. $ D_2 $ consists of
a finite sum of products of terms involving either a positive power of
$ x $ or a derivative of order at most $ k_0
- 1 $ of $y$ with respect to $ x_j $. In the first case we pick up an
$ \epsilon $ from
\page{67}
the powers of $x$ and in the second case we pick up a $ \delta $ from
the induction hypothesis. So we only need to determine the magnitude
of each term in $ D_2 $. Since each term in the expression for $ \psi
(y,x) $ is $ \O(\lambda^{\ell_0}) $, we only need to understand how
the powers of $ \lambda $ decrease when we differentiate a product
involving $ x^u $ and 
$$ 
\left( \frac{\partial^k y}{\partial x_j^k}\right)^p, \quad 1 
\leq k \leq k_0 - 2, 
$$ 
with respect to $ x_j $. Differentiating $ x^u $ gives
\page{68}
$ x_j^{-1} x^u $, losing $ \lambda^{-\frac{\ell_0}{m_j}} $ and
differentiating 
$$
\left(\frac{\partial^k y}{\partial x_j^k}\right)^p \ {\rm gives} \
\left(\frac{\partial^k y}{\partial x_j^k}\right)^{p - 1} \, 
\frac{\partial^{k + 1}y}{\partial x_j^{k + 1}},
\ \ {\rm losing} \  
\lambda^{- (1 - k \frac{\ell_0}{m_j})} \lambda^{1 - (k + 1)
\frac{\ell_0}{m_j}} = \lambda^{-\frac{\ell_0}{m_j}}  
$$ 
by induction. 
Therefore each term in $ D_2 $ is 
$ \O (\delta \lambda^{\ell_0 - k_0 \frac{\ell_0}{m_j}})$ 
and this finishes the proof of \eqref{4.10} and
thus (1) of the lemma. 

The proof of (2) follows in the same way as the proof of (1). The only
difference is that differentiating the term $ b_j x^{\alpha_j}_j
y^\beta $, $ b_j \neq 0 $, contributes a term $ b_j \alpha_j {!}\,
y^{\beta_j} \sim \lambda^{\beta_j} $ and 
\page{69}
so 
$$
D_2 \ \sim \ \lambda^{\beta_j} + \O (\delta \lambda^{\ell_0 - \alpha_j
\frac{\ell_0}{m_j}}) \ \sim \ \lambda^{\ell_0 - \alpha_j
\frac{\ell_0}{m_j}}
$$
since $ \frac{\alpha_j}{m_j} + \frac{\beta_j}{\ell_0} = 1 $. This
shows (2).

The proof of (3) follows similarly. We use induction on the partial
ordering $ u = (u_1, \ldots, u_{n - 2}) \leq \beta = (\beta_1, \ldots,
\beta_{n - 2}) $ if and only if $ u_j \leq \beta_j $, $ 1 \leq j \leq n
- 2 $. (1) and (2) show (3) is true for all pure derivatives, $ \beta
= (0, \ldots, \beta_j, \ldots, 0) $. The arguments used in proving (1)
and (2) show that if (3) is true for all
\page{70}
$ u \lvertneqq \beta $, then differentiating \eqref{4.5} shows
$$
0 \ = \ D_1 \, \frac{\partial^\beta y}{\partial x^\beta} \ + \ D_2
$$
where $ D_1 \sim \lambda^{\ell_0 - 1} $ and 
$$ 
D_2 = \O
\left(\lambda^{\ell_0 - \ell_0 \sum^{n - 2}_{j = 1} \frac{\beta_j}{m_j}}\right),
$$
proving (3).

Finally to prove (4), we first note that (3) implies that it is enough
to show (4) for any power $ p_\beta $. Again we use induction on the
partial ordering $ \leq $, supposing (4) is true for all $ u
\lvertneqq \beta $. Rewriting \eqref{4.4}
\page{71}
expresses \eqref{4.5} as
\begin{equation}
\eqlabel{4.11}
\lambda^{\ell_0} = Ay^{\ell_0} + \sum\stacksub{u}{u \leq \beta} a_u
(y) x^u + \O (|x|^{|\beta| + 1})
\end{equation}
where $ a_0 (y) = \O (y^{\ell_0 + 1}), \quad a'_0(y) = \O(y^{\ell_0})
$ and the $ a_u $'s are smooth. Taking the $ \beta $-th derivative of
\eqref{4.11} gives
$$
0 = [A \ell_0 y^{\ell_0 - 1} (0, \lambda) + a^{\prime}_0 (y (0, \lambda))] \,
\frac{\partial^\beta y}{\partial x^\beta}\, (0,\lambda) + C(\lambda)
$$
where $ C(\lambda) $ is a finite sum of terms of the form
$$
a^{(s)} (y (0, \lambda)) \prod\stacksub{u}{u \lvertneqq \beta}
\left( \frac{\partial^u y}{\partial x^u} (0, \lambda)\right)^{q_u}
$$
for some non-negative integers $ q_u $.
\page{72}
Here $ a(y) $ is either a power of $y$ or one of the $ a_u $'s. Using
the fact that $ y(0, \lambda) \sim \lambda $ and the inductive
hypothesis, we see that $ C(\lambda) \sim \lambda^p $ for some $p$ or
$ C(\lambda) = \O(\lambda^N) $ for every $N$. Since
$$
\frac{\partial^\beta y}{\partial x^\beta} (0,\lambda) = -
\frac{C(\lambda)}{[A \ell_0 y^{\ell_0 - 1} (0, \lambda) + a^{\prime}_0 (y (0,
\lambda))]}
$$
and $ a^{\prime}_0 (y) = \O(\lambda^{\ell_0}) $, we have shown (4) and this
finishes the proof of the lemma.
\end{proof}
\page{73}

We now turn back to the proof of \eqref{4.1} where we are examining
the integral in \eqref{4.9}. Let us write
$$
y (x , \lambda) = M_1 (x , \lambda) + M_2 (x , \lambda)
$$
where $ M_1(x , \lambda) $ is a polynomial in $ x_1 $ of degree $
\alpha_1 - 1 $ and $ M_2 $ is that part of the Taylor expansion of $
y(x , \lambda) $ in the variable $ x_1 $ that is $
\O(|x_1|^{\alpha_1}) $. We wish to replace the integral in \eqref{4.9}
by a similar integral where $ y (x , \lambda) $ is replaced by $
M_1(x , \lambda) $ and the $ \lambda $ integral is restricted
\page{74}
to 
$$
\lambda \leq \left(\frac{1}{|\eta|}\right)^{\frac{1}{\alpha_1 +
\kappa(\alpha_1)}} 
$$
where $ \kappa(\alpha_1) = 1 - \alpha_1 \frac{\ell_0}{m_1} $. Note that $
\alpha_1 + \kappa(\alpha_1) >  1 $.
Since
$$
\frac{\partial^{\alpha_1}y}{\partial x_1^{\alpha_1}} \sim \lambda^{1 -
\alpha_1 \frac{\ell_0}{m_1}}
$$
by part (2) of Lemma 1 and since $ \alpha_1 \geq 2 $, an application of
Van~der~Corput's lemma together with integration by parts shows
\page{75}
$$
\eqalign{
&\Biggl|\, \int_{\lambda \geq
\left(\frac{1}{|\eta|}\right)^{\frac{1}{\alpha_1 + \kappa(\alpha_1)}}}
e^{i \gamma \varphi(\lambda^{\ell_0})}
\int_{|x| \leq y^{\frac{\ell_0}{m} + \sigma}}
e^{i \xi \cdot x} e^{i \eta y (x , \lambda)} K(x_1 \lambda) \, dx \, d
\lambda \Biggr| \cr
&\qquad\leq C \left(\frac{1}{|\eta|}\right)^{1/\alpha_1}
\int_{\lambda  \geq  \left(\frac{1}{|\eta|}\right)^{{\frac{1}{\alpha_1 +
\kappa(\alpha_1)}}}} \left(\frac{1}{\lambda}\right)^{\frac{1}{\alpha_1}\left[1 -
\alpha_1 \frac{\ell_0}{m_1}\right]}
\int_{\R^{n - 2}} \frac{1}{|x|^n + \lambda^n} \, dx \, d \lambda \cr
%
%
&\qquad\leq C \left(\frac{1}{|\eta|}\right)^{\frac{1}{\alpha_1}} 
\int_{\lambda \geq
\left(\frac{1}{|\eta|}\right)^{\frac{1}{\alpha_1 + \kappa(\alpha_1)}}}
\frac{\lambda^{\frac{\ell_0}{m_1}}}{\lambda^{2 + \frac{1}{\alpha_1}}} \ d \lambda \  
%
\leq \ C \left(\frac{1}{|\eta|}\right)^{\frac{1}{\alpha_1}}
\, |\eta|^{\frac{1}{\alpha_1 + \kappa(\alpha_1)}\left[1 + \frac{1}{\alpha_1} -
\frac{\ell_0}{m_1}\right]} \cr
%
&\qquad= C \left(\frac{1}{|\eta|}\right)^{\frac{1}{\alpha_1}}
\, |\eta|^{\frac{1}{\alpha_1} \left[ \frac{\alpha_1 +
\kappa(\alpha_1)}{\alpha_1 + \kappa(\alpha_1)}\right]} \ = \ C.
}
$$

In the region 
$$
\lambda \leq \left(\frac{1}{|\eta|}\right)^{\frac{1}{\alpha_1 +
\kappa(\alpha_1)}}
$$
we would like to replace $ e^{i \eta y (x , \lambda)} $ by $ e^{i \eta
M_1 (x , \lambda)} $. We expect to be able to replace $ e^{i \eta y(x,
\lambda)} $ by $ e^{i \eta M_1 (x , \lambda)} $  with a bounded error
when 
$$ |x| \leq \left(\frac{1}{|\eta|
\lambda^{\kappa(\alpha_1)}}\right)^{\frac{1}{\alpha_1}} 
$$ 
since
\page{76}
$$| e^{i \eta y(x , \lambda)} - e^{i \eta M_1 (x , \lambda)}| \ \leq \ C
|\eta| \lambda^{1 - \alpha_1 \frac{\ell_0}{m_1}} |x|^{\alpha_1} \ \leq \ C
$$ 
when  
$
|x| \leq \left(\frac{1}{|\eta|
\lambda^{\kappa(\alpha_1)}}\right)^{\frac{1}{\alpha_1}}.
$
However in the complementary region, when
$$
\lambda \leq \left(\frac{1}{|\eta|}\right)^{\frac{1}{\alpha_1 +
\kappa(\alpha_2)}} 
\quad \text{ and } \quad |x| \geq \left(\frac{1}{|\eta| \lambda^{\kappa
(\alpha_1)}}\right)^{\frac{1}{\alpha_1}},
$$
$K$ is uniformly integrable. In fact
\page{77}
$$
\eqalign{
&\int_{\lambda \leq \left(\frac{1}{|\eta|}\right)^{\frac{1}{\alpha_1 +
\kappa (\alpha_1)}}} 
\int_{|x| \geq \left(\frac{1}{|\eta|
\lambda^{\kappa(\alpha_1)}}\right)^{\frac{1}{\alpha_1}}}
|K(x , \lambda)| \, d x \, d \lambda \cr
& \qquad \leq C \int_{\lambda \leq
\left(\frac{1}{|\eta|}\right)^{\frac{1}{\alpha_1 + \kappa(\alpha_1)}}} \,
\int_{|x| \geq \left(\frac{1}{|\eta|
\lambda^{\kappa(\alpha_1)}}\right)^{\frac{1}{\alpha_1}}}
\, \frac{1}{|x|^{n - 1}} \, d x \, d \lambda \cr
&\qquad \leq \int_{\lambda \leq
\left(\frac{1}{|\eta|}\right)^{\frac{1}{\alpha_1 + \kappa(\alpha_1)}}}
 ( |\eta| \lambda^{\kappa (\alpha_1)})^{\frac{1}{\alpha_1}} \, d
\lambda \cr
&\qquad  = C \int_{|\eta|^{\frac{1}{\alpha_1 + \kappa(\alpha_1)}} \lambda \leq 1}
\left(|\eta|^{\frac{1}{\alpha_1 + \kappa(\alpha_1)}} \lambda
\right)^{\frac{\kappa(\alpha_1)}{\alpha_1} + 1} \, \frac{d
\lambda}{\lambda} \ \leq \ C}
$$
since
$$
\frac{\kappa(\alpha_1)}{\alpha_1} = \frac{1 - \alpha_1
\frac{\ell_0}{m_1}}{\alpha_1} >  -1.
$$
Replacing $ e^{i \eta y(x , \lambda)} $ by $ e^{i \eta M_1(x ,
\lambda)} $ when 
$$ 
\lambda \leq
\left(\frac{1}{|\eta|}\right)^{\frac{1}{\alpha_1 + \kappa(\alpha_1)}} 
\quad \text{ and } \quad
|x| \leq \left(\frac{1}{|\eta|
\lambda^{\kappa(\alpha_1)}}\right)^{\frac{1}{\alpha_1}} 
$$ 
creates an error at most 
\page{78}
$$
\eqalign{
&C |\eta| \int_{\lambda \leq
\left(\frac{1}{|\eta|}\right)^{\frac{1}{\alpha_1 + \kappa(\alpha_1)}}} \,
\lambda^{\kappa(\alpha_1)}
\int_{|x| \leq \left(\frac{1}{|\eta|
\lambda^{\kappa(\alpha_1)}}\right)^{\frac{1}{\alpha_1}}}
\, \frac{|x|^{\alpha_1}}{|x|^{n - 1}} \, d x \, d \lambda \cr
& \qquad \leq  C |\eta| \int_{\lambda \leq
\left(\frac{1}{|\eta|}\right)^{\frac{1}{\alpha_1} + \kappa(\alpha_1)}} \,
\frac{\lambda^{\kappa(\alpha_1)}}{\left(|\eta|
\lambda^{\kappa(\alpha_1)}\right)^{\frac{1}{\alpha_1} (\alpha_1 - 1)}} \, d
\lambda \cr
%
& \qquad \leq C \int_{|\eta|^{\frac{1}{\alpha_1 + \kappa(\alpha_1)}} \lambda
\leq 1}%
\, \left(|\eta|^{\frac{1}{\alpha_1 + \kappa(\alpha_1)}} \lambda
\right)^{\frac{\alpha_1 + \kappa(\alpha_1)}{\alpha_1}} \, \frac{d
\lambda}{\lambda} \ \leq \ C
}
$$
since 
$$
\frac{\alpha_1 + \kappa(\alpha_1)}{\alpha_1} = 1 - \frac{\ell_0}{m_1} +
\frac{1}{\alpha_1} >  0.
$$
Therefore
\page{79}
$$
I = \int\stacksub{0 \leq \lambda \leq 1}{\lambda \leq 
\left(\frac{1}{|\eta|}\right)^{\frac{1}{\alpha_1 + K(\alpha_1)}}}
e^{i \gamma \varphi(\lambda^{\ell_0})}
\int_{|x| \leq y^{\frac{\ell_0}{m} + \sigma}}
e^{i \xi \cdot x} e^{i \eta M_1(x , \lambda)} \, K(x , \lambda) \, dx
\, d \lambda \ + \ \O(1).
$$
Since
$$
M_1(x , \lambda) = \sum^{\alpha_1 - 1}_{k = 0} \, \frac{1}{k{!}} \,
\frac{\partial^k y}{\partial x_1^k} (0, x_2, \ldots, x_{n - 2})
x_1^k,
$$
we see that for $ 2 \leq j \leq n - 2 $,
$$
\frac{\partial^{\alpha_j}M_1}{\partial x_j^{\alpha_j}} \, (x ,
\lambda) = \frac{\partial^{\alpha_j}y}{\partial x_j^{\alpha_j}} (0,
x_2, \ldots, x_{n - 2}) + \sum^{\alpha_1 - 1}_{k = 1} \, \frac{1}{k{!}}\,
\frac{\partial^{k + \alpha_j}y}{\partial x_1^k \partial x_j^{\alpha_j}} (0,x_2,
\ldots, x_{n - 2}) x_1^k.
$$
Also since $ |x_1| \leq \epsilon \lambda^{\frac{\ell_0}{m_1}} $, we have by
part (3) of Lemma \ref{1},
$$
\left| \sum^{\alpha_1 - 1}_{k = 1} \, \frac{1}{k{!}}\, \frac{\partial^{k +
\alpha_j} y}{\partial x_1^k \partial x_j^{\alpha_j}} (0, x_2, \ldots,
x_{n - 2}) x_1^k \right| \leq \epsilon \lambda^{k \frac{\ell_0}{m_1}}
\lambda^{ 1 - \ell_0 \left( \frac{k}{m_1} +
\frac{\alpha_j}{m_j}\right)}
\leq \epsilon \lambda^{1 - \ell_0 \frac{\alpha_j}{m_j}},
$$
\page{80}
and since 
$$
\frac{\partial^{\alpha_j} y}{\partial x_j^{\alpha_j}} \ \sim \ \lambda^{ 1
- \ell_0 \frac{\alpha_j}{m_j}}
$$
by part (2) of Lemma \ref{1}, we conclude that for $ 2 \leq  j \leq n - 2 $, 
$$
\frac{\partial^{\alpha_j} M_1}{\partial x_j^{\alpha_j}} (x ,  \lambda)
\ \sim \ \lambda^{ 1 - \ell_0 \frac{\alpha_j}{m_j}}.
$$
Therefore we may proceed in the same manner to find that up to a
bounded error
$$
I = \int_{\lambda \leq
\frac{1}{|\eta|^\delta}}
e^{i \gamma  \varphi(\lambda^{\ell_0})} \int_{|x| \leq
\lambda^{\frac{\ell_0}{m} + \sigma}}
e^{i \xi \cdot x} e^{i \eta Q(x , \lambda)} \, K(x , \lambda) \, dx \,
d \lambda
$$
\page{81}
for some $ 0 < \delta < 1 $ where
\marginpar{Punctuation?}
$$
Q(x , \lambda) = \sum\stacksub{\beta}{\beta_j \leq \alpha_j - 1} \,
\frac{1}{\beta{!}} \, \frac{\partial^\beta y}{\partial x^\beta} \,
(0 , \lambda) \, x^\beta.
$$
By part (4) of Lemma \ref{1}, we have for each $ \beta $ with $
\beta_j \leq \alpha_j - 1 $, $ 1 \leq j \leq n - 2 $, either 
$$
\frac{\partial^\beta y}{\partial x^\beta} (0 , \lambda) \ \sim \
\lambda^{p_\beta} 
$$
for some 
$ p_\beta >  - |\beta| $ 
or $$ 
\frac{\partial^\beta y}{\partial x^\beta}\, (0, \lambda) \ = \ \O \, (\lambda^N) $$ 
for every $N$. If the latter
occurs, then up to a bounded error, we may clearly replace $ e^{i \eta
\frac{\partial^\beta y}{\partial x^\beta} (0 , \lambda) x^\beta} $ by
$ 1 $. When the former occurs, that is, when $ \frac{\partial^\beta
y}{\partial x^\beta} (0 , \lambda) $ behaves like a power of $ \lambda
$, we can repeat the above argument to see that for some (other) $
\delta $, $ 0 < \delta < 1 $,
\page{82}
$$
I = \int_{\lambda \leq \frac{1}{|\eta|^\delta}}
e^{i \gamma \varphi(\lambda^{\ell_0})}
e^{i \eta y (0 , \lambda)}
\int_{|x| \leq \lambda^{\frac{\ell_0}{m} + \sigma}}
e^{i \sum\limits^{n - 2}_{j = 1} (\xi_j + \lambda^{p_j} \eta_j) x_j}
\, K (x , \lambda) \, dx \, d \lambda \ + \ \O (1).
$$
For each $ x_j $ integral, by splitting the $ \lambda $ integral where
$ \lambda $ is smaller or larger than 
$ | \xi_j / \eta_j|^{\frac{1}{p_j}} $, 
we can once again repeat the same argument to conclude that
\marginpar{Ms. illegible at ??}
$$
\eqalign{
I &= \int_{\lambda \leq \min(\frac{1}{|\xi|}, \frac{1}{|\eta|^{b}})}
e^{i \gamma \varphi(\lambda^{\ell_0})}
e^{i \eta y (0 , \lambda)}
\int_{|x| \leq \lambda^{\frac{\ell_0}{m} + \sigma}}
\, K(x , \lambda) \, dx \, d \lambda \ + \ \O(1)\cr
& = \int_{\lambda \leq \min (\frac{1}{|\xi|}, \frac{1}{|\eta|^b})}
e^{i \gamma \varphi(\lambda^{\ell_0})}
e^{i \eta y (0 , \lambda)}
\, \frac{1}{\lambda} \int_{\R^{n - 2}} K(z, 1 ) \, dz \, d \lambda \ +
\ \O(1)
}
$$
for \page{83} some $ 0 < b < 1 $. Thus the proof 
of \eqref{4.1} will be finished once we establish the identity
\begin{equation}
\eqlabel{4.12}
\int_{\R^{n - 2}} K(x , 1) \, dx = \int_{\Sigma^+} K(\omega) \, d
\sigma(\omega).
\end{equation}
This is done by making the change of variables
$$
x_j = \frac{s_j}{1 - |s|^2}, \quad 1 \leq j \leq n - 2.
$$
\page{84}
In evaluating the Jacobian of this change of variables, we need to
observe that if an $ r \times r $ matrix $ (\alpha_{j,k}) $ is defined
by $ \alpha_{j,k} = s_j s_k $ for $ j \neq k $ and $ \alpha_{j,j} = 1
- |s|^2 - s^2_j $, then
$$
\det (\alpha_{j,k}) = (1 - |s|^2)^{r - 1}.
$$
This calculation was shown to us by A.~Carbery and
is carried out in the appendix. This establishes
\eqref{4.12} and finishes the proof of \eqref{4.1}. The proof of
\eqref{4.2} is similar. It remains to prove \eqref{4.3}.
\page{85}

Suppose first that there is no constant $ C_0 $ so that
$\bar{\varphi}' (C_0 \lambda) \geq 2 \bar{\varphi}' (\lambda) $ for
$ 0 < \lambda \leq 1 $. Then there exists a sequence of points $
\lambda_j \searrow 0 $ such that
$$
\frac{\lambda_j \bar{\varphi}' (\lambda_j)}{\lambda_j \bar{\varphi}'
(\lambda_j) - \bar{\varphi}(\lambda_j)} \rightarrow \infty.
$$
See, e.g., \cite{NVWW}. Let
$$
\gamma_j = \frac{\pi}{4} \, \frac{1}{\lambda_j \bar{\varphi}^1
(\lambda_j) - \bar{\varphi}(\lambda_j)}, \quad \eta_j = \gamma_j
\bar{\varphi}' (\lambda_j) ,
$$
and choose $ \xi_j $ so that
$$
\frac{1}{\xi_j} = \min \left( \lambda_j, \frac{1}{\eta_j^b},
\frac{1}{\eta_j^{1 + \epsilon}}\right)
$$
\page{86}
where $ \epsilon >  0 $ is chosen so that $ q (\lambda) = \lambda + \O
(\lambda^{ 1 + \epsilon}) $. Then
$$
\eqalign{
\left|\int_0^{\frac{1}{\xi}} e^{i \gamma_j \bar{\varphi} (\lambda)}
\sin(\eta_j q(\lambda))\, \frac{d \lambda}{\lambda} \right|
&= \left| \int^{\frac{1}{\xi}}_{\frac{1}{\eta_j}} e^{i (\gamma_j
\bar{\varphi}(\lambda) - \eta_j \lambda)} \, \frac{d
\lambda}{\lambda} \right| \ + \ \O(1)\cr
&\geq A \log \left(\frac{\eta_j}{\xi_j}\right) \geq A
\log(\lambda_j \eta_j) \rightarrow + \infty
}
$$
for some $ A >  0 $ since $ \lambda_j \eta_j \rightarrow \infty $ and
$$
0 \leq \eta_j \lambda - \gamma_j \bar{\varphi}(\lambda) \leq
\frac{\pi}{4}
$$
for all $ 0 \leq \lambda \leq \lambda_j $.
\page{87}

Finally let us turn to the proof of the sufficiency of \eqref{4.3} and
assume $ \bar{\varphi}(0) = \bar{\varphi}'(0) = 0 $, $
\bar{\varphi}' (C_0 t) \geq 2 \bar{\varphi}' (t) $
\marginpar{Are these o's or 0's?} for some $ C_0
\geq 1 $. It suffices to show that the integral 
\begin{equation}
\eqlabel{4.13}
\I = \int_{\frac{1}{\eta} \leq t \leq 1} e^{i \gamma \bar{\varphi}(t)}
e^{- i \eta q(t)} \, \frac{dt}{t}
\end{equation}
is uniformly bounded in $ \gamma, \eta >  0 $. First assume $ 10
\gamma >  \eta $. Choosing $ t_0 $ such that $ \bar{\varphi}' (t_0) =
\frac{\eta}{\gamma} $ we write
$$
\I = \int_{\frac{1}{\eta} \leq t \leq \frac{t_0}{C_0}} +
\int_{\frac{t_0}{C_0} \leq t \leq C_0 t_0} + \int_{C_0 t_0 \leq  t \leq 1}
= A + B + D.
$$
\page{88}
For $ \frac{1}{\eta} \leq t \leq \frac{t_0}{C_0} $, $ \frac{d}{dt}
(\eta t - \gamma \bar{\varphi} (t) )  = \eta - \gamma \bar{\varphi}'
(t) \geq \eta - \gamma \bar{\varphi}' (\frac{t_0}{C_0}) \geq
\frac{\eta}{2} $, and so integrating by parts shows
$$
\eqalign{
|A| &\leq  \frac{1}{\eta} \int_{\frac{1}{\eta} \leq t \leq 1} \eta |
q'(t) - 1| \, \frac{dt}{t} + \frac{1}{\eta} \int_{\frac{1}{\eta} \leq
t} \, \frac{1}{t^2} \, dt  \ + \ C \cr
& \leq C \int^1_0 t^\epsilon \, \frac{dt}{t} + C \leq C.
}
$$
Also
$$
|B| \leq \int_{\frac{t_0}{C_0} \leq t \leq C_0 t_0} \, \frac{1}{t} \,
dt \leq 2\log(C_0).
$$
For $ C_0 t_0 \leq t $, $ \frac{d}{dt} (\gamma \bar{\varphi}(t) - \eta
t) = \gamma \bar{\varphi}' (t) - \eta \geq \frac{\gamma}{2}
\bar{\varphi}' (t) $, 
\page{89}
and so integrating by parts show
$$
\eqalign{
|D| & \leq \frac{1}{\gamma} \int_{\frac{1}{\eta} \leq t_0 \leq t}
\left( \frac{\bar{\varphi}''(t)}{t [\bar{\varphi}'(t)]^2} +
\frac{1}{\bar{\varphi}' (t) t^2} + \frac{\eta | q'(t) - 1|
}{\bar{\varphi}' (t) t}\right) \, dt + \frac{\eta}{\gamma}\frac{1}{\bar{\varphi}'
(t_0 )} \cr
& \leq \frac{\eta}{\gamma} \int_{t_0 \leq t}
\frac{\bar{\varphi}''(t)}{[\bar{\varphi}^{\prime}(t)]^2} \, dt +
\frac{1}{\bar{\varphi}'(t_0)} \left[\frac{\eta}{\gamma} + C
\frac{\eta}{\gamma} \int^1_0 \, \frac{t^\epsilon}{t}\, dt \right] + \frac{\eta}{
\gamma}\frac{1}{\bar{\varphi}' (t_0 )} \cr
& \leq C \frac{\eta}{\gamma} \, \frac{1}{\bar{\varphi}' (t_o)} \leq C
}
$$
since $ \frac{\eta}{\gamma} = \bar{\varphi}'(t_0) $. Next suppose
$ 10 \gamma \leq \eta $. Then in a neighborhood of the origin,
\page{90}
$ \frac{d}{dt} [\eta t - \gamma \bar{\varphi}(t) ] \geq \frac{\eta}{2}
$, and so integrating by parts shows
$$
\eqalign{
| \I | &\leq \frac{1}{\eta} \left[ \int_{\frac{1}{\eta} \leq t} \eta |
q' (t) - 1| \, \frac{dt}{t} + \int_{\frac{1}{\eta} \leq t} \,
\frac{1}{t^2}\right] + C \cr
& \leq C \int^1_0 t^\epsilon \frac{dt}{t} + C \leq C.
}
$$
This completes the proof of Theorems \ref{2} and \ref{3}.
\qed

\page{91}

\section{Proof of Theorem \ref{4}}
We will prove the $ L^p $ boundedness of the maximal function
$$
Mf(x', x_n) = \sup_{0 < h \leq 1} \frac{1}{h^{n - 1}} \left| \int_{|t|
\leq h} f(x ' - t, x_n - \varphi(\psi(t)) \, dt \right|.
$$
The proof for the singular integral is similar.

When $ E_{\ell_0} = \{ 0 \} $, the main term $ P(t) $, in the
decomposition \eqref{3.1} for $ \psi(t) $, $ \psi(t) = P(t) + R(t) $,
is a positive homogeneous polynomial of degree $ \ell_0 $. $ R(t) $
consists of
\page{92}
all the terms in the Taylor expansion of $ \psi $ with degree greater
than $ \ell_0 $. The proof of $ L^p $ boundedness for $M$ in the case
$ P(t) = |t|^2 $ and $ R(t) \equiv 0 $ is carried out in  \cite{KWWZ}.
We will see that slight modifications of the arguments given in
\cite{KWWZ} work for the general case.

It will be convenient for us to use polar coordinates with respect to
the surface $ P(\omega) = 1 $. That is, every $ t \neq 0 \in \R^{n -
1} $ can be written uniquely
\page{93}
as $ t = r \omega $ where $ r >  0 $ and $ P(\omega) = 1 $. We also
introduce a norm $ \|\cdot \| $ so that $ \|t\| = \|rw \| = r $. Since
the Euclidean norm of $ \omega $, $ |\omega | $, is bounded above and
below as $ \omega $ runs over the surface $ P(\omega) = 1 $, it is
clear that the maximal function $ Mf(x) $ is pointwise comparable to
the maximal function defined in terms of averages with respect to the
norm $ \| \cdot \| $. Therefore it suffices to consider
$$
\cal{M} f(x) = \sup_{k >  0} 2^{k(n - 1)} \left| \int \chi (2^k \|t\|)
f(x - \Gamma(t)) \, dt \right|
\stackrel{\text{def}}{=} \sup_{k >  0} |f * d \mu_{k}
(x) |
$$
\marginpar{Is this really a $*$?}
\page{94}
where $ \chi $ is a smooth cut-off function supported in $ [1,2] $ and
chosen so that
$$
2^{k(n - 1)} \int_{\R^{n - 1}} \chi (2^k \|t\|) \, dt \equiv 1 .
$$

To prove $ L^p $ bounds for $\cal{M} $ we introduce dilations $ \{
\delta(t) \}_{t >  0} $, defined by $ \delta(t) (\xi, \gamma) = ( t
\xi, \bar{\varphi}(t) \gamma ) $. Although the ``balls'' generated
with respect to these dilations do not in general form
\page{95}
a space of homogeneous type with respect to Lebesgue measure, an
appropriate singular integral and Littlewood--Paley theory for the
dilations $ \{ \delta(t) \}_{t >  0} $ has been worked out in
\cite{CCVWW}. Using this theory and well-known techniques, following
the arguments detailed in \cite{KWWZ}, we reduce ourselves to proving
two basic estimates for the Fourier transform of the measures $ \{
d \mu_k \} $ defined above:
\page{96}
\begin{equation}
\eqlabel{5.1}
|\widehat{d \mu_k} (\xi, \gamma) -1 | \leq C |\delta (2^{- k + 3})
(\xi, \gamma) |,
\end{equation}
and
\begin{equation}
\eqlabel{5.2}
| \widehat{d \mu_k} (\xi, \gamma)| \leq C | \delta (2^{-k - 1})(\xi,
\gamma) |^{- \epsilon}
\end{equation}
for some $ \epsilon >  0 $. Using polar coordinates $ t = r \omega $,
\begin{equation}
\eqlabel{5.3}
\widehat{d \mu_k} (\xi, \gamma) = 2^{k(n - 1)}
\int_{\R}\!\int_{P(\omega) = 1} \chi (2^k r) e^{i \xi r \cdot \omega}
e^{i \gamma \varphi (\psi(r \omega))} r^{n - 2} h(\omega) \, d \omega
\, dr
\end{equation}
where $ h (\omega) $ is some smooth function. Since $ \psi(r \omega) =
r^{\ell_0} + \O(r^{\ell_0 + 1}) $, we have for $ k >  0 $ large, $
\varphi (\psi (r \omega)) \leq \bar{\varphi} (2^{-k + 3}) $ when $
2^{-k} \leq r \leq 2^{-k + 1} $. Therefore
\page{97} 
$$
| \widehat{d \mu_k} (\xi, \gamma) - 1 | \leq C [ 2^{-k} |\xi| +
\bar{\varphi} (2^{-k + 3})|\gamma| ] \leq C | \delta (2^{-k + 3})(\xi,
\gamma)|,
$$
establishing \eqref{5.1}. To prove \eqref{5.2} we make the change of
variables
\begin{equation}
\eqlabel{5.4}
\lambda^{\ell_0} = \psi(r \omega) =  r^{\ell_0} + R( r \omega) 
\end{equation}
in the $r$ integral in \eqref{5.3} for fixed $ \omega $. For $ k >  0
$ large this is a good change of variables and so
\begin{equation}
\eqlabel{5.5}
\widehat{d \mu_k} (\xi, \gamma) = 2^{k(n - 1)} \int_\R e^{i \gamma
\bar{\varphi}(\lambda)} \int_{P(\omega) = 1} e^{i r(\lambda, \omega)
\xi \cdot \omega}
\chi (2^k r (\lambda, \omega)) \, \frac{\partial r}{\partial \lambda}
\,
r^{n - 2}{(\lambda, \omega)}\, h(\omega) \, d \omega \, d \lambda.
\end{equation}
From \eqref{5.4} one easily deduces the following estimates on the
derivatives of $ r(\lambda,  \omega) $:
\begin{equation}
\eqlabel{5.6}
r(\lambda, \omega) = \lambda + \O(\lambda^2), \quad \frac{\partial
r}{\partial \lambda} = 1 + \O(\lambda),
\end{equation}
\begin{equation}
\eqlabel{5.7}
\nabla_\omega r = \O(\lambda), \quad \frac{\partial r}{\partial
\lambda \partial \omega} = \O(1),
\end{equation}
\begin{equation}
\eqlabel{5.8}
\frac{\partial^2 r}{\partial \lambda^2} = \O \left(\frac{1}{\lambda}\right).
\end{equation}
Since $ P(\omega) = 1 $ is of finite type, we can argue as in section
3 to find an $ \epsilon >  0 $ such that 
$$
\left| \int_{P(\omega) = 1} e^{i r(\lambda, \omega) \xi \cdot \omega}
h(\omega) \, d \omega \right| \leq C \frac{1}{|\lambda \xi|^{2
\epsilon}}.
$$
Now integrating by parts, using \eqref{5.6} and \eqref{5.7}, shows
\page{99}
$$
|\widehat{d \mu_k} (\xi, \gamma) | \leq C \left(\frac{1}{2^{-k - 1}
|\xi|}\right)^{2 \epsilon},
$$
establishing \eqref{5.2} if 
$$ 
\sqrt{|\gamma| \bar{\varphi}(2^{- k - 1})} \leq C 2^{- k - 1} |\xi|.
$$
On the other hand, if 
$$
C 2^{- k - 1} |\xi| \leq \sqrt{|\gamma| \bar{\varphi}(2^{- k - 1})},
$$
we perform the $ \lambda $ integration first, writing \eqref{5.5} as
$$
\widehat{d \mu_k}(\xi, \gamma) = 2^{k(n - 1)} \int_{P(\omega) = 1}
h(\omega) \int_\R e^{i[\gamma \bar{\varphi}(\lambda) + \lambda \xi \cdot
\omega]} e^{i[r(\lambda, w) - \lambda] \xi \cdot \omega }
\chi (2^k r) \, \frac{\partial r}{\partial \lambda} \, r^{n - 2} \, d
\lambda \, d \omega.
$$
For $ 2^{- k} \leq r (\lambda, \omega) \leq 2^{ - k + 1} $, we have
$$
\left| \frac{\partial}{\partial \lambda} [\gamma
\bar{\varphi}(\lambda) + \lambda \xi \cdot \omega ] \right| \geq
\frac{|\gamma|}{2} \bar{\varphi}' (2^{- k - 1}) \geq |\gamma| \,
\frac{\bar{\varphi} (2^{- k - 1})}{2^{- k}}
$$
since $ 2^{- k} |\xi| \ll |\gamma| \bar{\varphi} (2^{- k - 1}) $. Thus
\page{100}
integrating by parts, using \eqref{5.6} and \eqref{5.8}, shows
$$
\eqalign{
|\widehat{d \mu_k} (\xi, \gamma)| &\leq C \left[ \frac{1}{|\gamma|
\bar{\varphi} (2^{- k - 1})} + \frac{|\xi| 2^{ - k}}{|\gamma|
\bar{\varphi} (2^{ - k - 1})} \right]\cr
&\leq C \frac{1}{\sqrt{|\gamma| \bar{\varphi} (2^{- k - 1})}} \cr
&\leq C
\frac{1}{\sqrt{|\delta (2^{- k - 1})(\xi, \gamma)|}}
}
$$
since $ C 2^{- k}|\xi| \leq \sqrt{|\gamma| \bar{\varphi} (2^{- k -
1})} $. 

This completes the proof of \eqref{5.1} and \eqref{5.2} from which the
$ L^p $ boundedness of the maximal function follows as in \cite{KWWZ}.
\page{101}
\section{Proof of Theorem 5}
We need  to show that the multiplier for $H$,
\begin{equation}
\eqlabel{6.1}
m(\xi, \gamma) = \int\!\int\stacksub{|t| \leq 1}{t \in \R^2} e^{i
\gamma \varphi(\psi(t))} e^{i \xi \cdot t} K(t) \, dt
\end{equation}
is uniformly bounded for $ \xi \in \R^2 $ and $ \gamma \in \R $.
Introducing polar coordinates with respect to the convex curve $
\psi(t) = 1 $, we may write \eqref{6.1} as
\begin{equation}
\eqlabel{6.2}
\int^1_0 e^{i \gamma \varphi(r)} \, \frac{1}{r} \int_{\psi(\omega) =
1} e^{i r \xi \cdot \omega} K(\omega)\, h(\omega) \, d \omega \, dr
\end{equation}
for some smooth function $ h(\omega) $. The
\page{102}
argument used in the proof of Theorem \ref{1} to establish
\eqref{3.10} shows
$$
\int_{\psi(\omega) = 1} K(\omega) \, h(\omega) \, d \omega = 0
$$
and so the part of the integral in \eqref{6.2} where $ r \leq
\frac{1}{|\xi|} $ is at most
$$
C \int_{r \leq \frac{1}{|\xi|}} \, \frac{1}{r}
\, \left|\int_{\psi(\omega) = 1} (e^{i r \xi \cdot \omega} - 1) K(\omega)
h(\omega) \, d \omega \right| \, dr \leq C |\xi| \int_{r \leq
\frac{1}{|\xi|}} \, dr \leq C.
$$

For the region where $ r \geq \frac{1}{|\xi|} $, we observe that the
inner integral
\page{103}
in \eqref{6.2} is the Fourier transform of a smooth density on the
convex curve $ \psi (\omega) = 1 $ evaluated at $ r \xi $. This Fourier
transform can be estimated in terms of the ``balls,'' $ E(t,
\epsilon)$, introduced in section 1. In fact
$$
\left|\int_{\psi  (\omega) = 1} e^{ i r \xi \cdot \omega} \, K(\omega)
h(\omega) \, d \omega \right| \leq C \left[ | E (t_1 (\xi),
\frac{1}{r|\xi|})| + | E (t_2(\xi), \frac{1}{r|\xi|})| \right]
$$
where $ t_1(\xi) $ and $ t_2(\xi ) $ are the two points on the curve $
\psi (t) = 1 $ whose tangent lines are normal to $ \xi $. See
\cite{BNW}. Therefore the
\page{104}
part of the integral in \eqref{6.2} where $ r \geq \frac{1}{|xi|} $ can
be estimated by
$$
C \sup\stacksub{t}{\psi(t) = 1} \int_{\frac{1}{|\xi|} \leq r} | E( t,
\frac{1}{r|\xi|})| \, \frac{dr}{r} \leq 
C \sup\stacksub{t}{\psi(t) = 1} \int_{0}^{1} |E (t, \delta)| \,
\frac{d \delta}{\delta}.
$$
Hence the multiplier $ m (\xi, \gamma) $ is uniformly bounded in $
\xi $ and $ \gamma $ if the quantity
$$
\sup\stacksub{t}{\psi(t) = 1} \int^1_0 |E (t, \delta)| \,
\frac{d \delta}{\delta}
$$
is finite. This completes the proof of Theorem \ref{5}.
\qed
\section{Appendix}
In this appendix we will compute the determinant of
an $r\times r$ matrix $A = \{ \alpha_{j,k} \}$ of the form $A = cI + B$
where $B = \{ b_{j,k} \}$ and $b_{j,k} = b\, s_j t_k$. We will
show that 
\begin{equation}
\eqlabel{7.1}
\det(A) \ = \ c^r \ + \ c^{r-1}\, b \, \sum\limits_{j=1}^{r} s_j t_j .
\end{equation}
For the example we need in this paper, $\alpha_{j,k} = s_j s_k$ for
$j\ne k$ and $\alpha_{j,j} = 1 - |s|^2 + s_j^2$. Therefore taking
$t_j = s_j$, $c = 1 - |s|^2$ and $b = 1$ in the above formula (7.1) gives
us the desired result $\det(A) = (1 - |s|^2 )^{r-1}$ in this case.
To prove (7.1) first note that as a function of $s=(s_1,\ldots ,s_r )$,
$\det(A)$ is an affine function in each of the variables $s_j$ separately.
Also computing any pure mixed derivative, e.g., $\frac{\partial^3}{\partial s_1
\partial s_2 \partial s_3 }$, of $\det(A)$ gives rise to two or more
rows being identical and therefore zero. Hence expanding $\det(A)$
in its Taylor series in $s$ about the origin, we see that (7.1) follows
from the fact that for each $1\le j \le r$, the partial derivative
of $\det(A)$ with respect to $s_j$ at the origin is $c^{r-1}\, b\, t_j$.
This is a straightforward computation.
\qed

\marginpar{Need title for \cite{CZ} reference.}

\bibliographystyle{alpha}

\begin{thebibliography}{NVWW83}

\bibitem[BNW]{BNW}
J.~Bruna, A.~Nagel, and S.~Wainger.
\newblock Convex hypersurfaces and {F}ourier transforms.
\newblock {\em Ann. of Math.}, 127:333--365, 1988.

\bibitem[CCVWW]{CCVWW}
Anthony Carbery, Michael Christ, James Vance, Stephen Wainger, and David~K.
  Watson.
\newblock Operators associated to flat plane curves: ${L}\sp p$ estimates via
  dilation methods.
\newblock {\em Duke Math. J.}, 59(3):675--700, 1989.

\bibitem[KWWZ]{KWWZ}
Weon-Ju Kim, Stephen Wainger, James Wright, and Sarah Ziesler.
\newblock Singular integrals and maximal functions associated to surfaces of
  revolution.
\newblock {\em Bull. London Math. Soc.}, 28(3):291--296, 1996.

\bibitem[NVWW]{NVWW}
A.~Nagel, J.~Vance, S.~Wainger, and D.~Weinberg.
\newblock Hilbert transforms for convex curves.
\newblock {\em Duke Math. J.}, 50:735--744, 1983.

\bibitem[Sch]{Sc}
Helmut Schulz.
\newblock Convex hypersurfaces of finite type and the asymptotics of their
  {F}ourier transforms.
\newblock {\em Indiana Univ. Math. J.}, 40(4):1267--1275, 1991.

\bibitem[SWWZ]{SWWZ}
A.~Seeger, S.~Wainger, J.~Wright, and S.~Ziesler.
\newblock Classes of singular integrals along curves and surfaces.
\newblock {\em {T}rans. {A}mer. {M}ath. {S}oc.}
\newblock To appear.

\bibitem[St]{St}
E.~Stein.
\newblock {\em Harmonic Analysis}
\newblock Princeton University Press, Princeton New Jersey, 1993.

\end{thebibliography}

\newcommand{\etalchar}[1]{$^{#1}$}

\end{document}